\documentclass{amsart}
\usepackage[all]{xy}
\usepackage{verbatim}
\usepackage{color}
\usepackage{amsthm}
\usepackage{amssymb}
\usepackage[colorlinks=true]{hyperref}
\usepackage{graphicx}
\usepackage{mathtools}
\usepackage[utf8]{inputenc}

\usepackage{footmisc}
\usepackage{stmaryrd}

\numberwithin{equation}{section}

\newtheorem{theorem}[equation]{Theorem}
\newtheorem*{theorem*}{Theorem}

\newtheorem*{conjecture*}{Mamma Conjecture}
\newtheorem*{conjecture1*}{Mamma Conjecture (revisited)}

\newtheorem{corollary}[equation]{Corollary}
\newtheorem*{corollary*}{Corollary}

\theoremstyle{remark}

\newtheorem{example}[equation]{Example}

\newtheorem{notation}[equation]{Notation}

\theoremstyle{remark}
\newtheorem{remark}[equation]{Remark}

\newcommand{\cA}{{\mathcal A}}
\newcommand{\cB}{{\mathcal B}}
\newcommand{\cC}{{\mathcal C}}
\newcommand{\cD}{{\mathcal D}}
\newcommand{\cE}{{\mathcal E}}
\newcommand{\cF}{{\mathcal F}}

\newcommand{\cL}{{\mathcal L}}
\newcommand{\cM}{{\mathcal M}}
\newcommand{\cN}{{\mathcal N}}
\newcommand{\cO}{{\mathcal O}}
\newcommand{\cP}{{\mathcal P}}

\newcommand{\cT}{{\mathcal T}}

\newcommand{\cV}{{\mathcal V}}

\newcommand{\cX}{{\mathcal X}}

\newcommand{\cZ}{{\mathcal Z}}

\newcommand{\bbA}{\mathbb{A}}

\newcommand{\bbC}{\mathbb{C}}

\newcommand{\bbF}{\mathbb{F}}
\newcommand{\bbG}{\mathbb{G}}

\newcommand{\bbL}{\mathbb{L}}
\newcommand{\bbN}{\mathbb{N}}
\newcommand{\bbP}{\mathbb{P}}
\newcommand{\bbR}{\mathbb{R}}

\newcommand{\bbQ}{\mathbb{Q}}
\newcommand{\bbZ}{\mathbb{Z}}

\DeclareMathOperator{\id}{id}

\DeclareMathOperator{\NChow}{NChow} 
\DeclareMathOperator{\NNum}{NNum} 

\DeclareMathOperator{\Num}{Num} 

\newcommand{\dgcat}{\mathrm{dgcat}} 

\newcommand{\perf}{\mathrm{perf}}

\newcommand{\Chow}{\mathrm{Chow}}

\newcommand{\dg}{\mathrm{dg}}

\newcommand{\uHom}{\underline{\mathrm{Hom}}}
\newcommand{\Hom}{\mathrm{Hom}}

\newcommand{\dgHo}{\mathrm{H}^0}

\newcommand{\Hmo}{\mathrm{Hmo}}
\newcommand{\op}{\mathrm{op}}

\newcommand{\too}{\longrightarrow}

\newcommand{\ie}{\textsl{i.e.}\ }
\newcommand{\eg}{\textsl{e.g.}}

\newcommand{\loc}{\mathrm{loc}}
\newcommand{\Nmot}{\mathrm{NMot}}
\newcommand{\Nmix}{\mathrm{NMix}}
\newcommand{\NHom}{\mathrm{NHom}}
\newcommand{\NVoev}{\mathrm{NVoev}}

\newcommand{\vphi}{{\varphi\!}{/}{\!\sim}}

\let\oldmarginpar\marginpar
\def\marginpar#1{\oldmarginpar{\tiny #1}}

\def\multiset#1#2{\ensuremath{\left(\kern-.3em\left(\genfrac{}{}{0pt}{}{#1}{#2}\right)\kern-.3em\right)}}

\begin{document}

\title[Recent developments on noncommutative motives]{Recent developments on noncommutative motives}
\author{Gon{\c c}alo~Tabuada}

\address{Gon{\c c}alo Tabuada, Department of Mathematics, MIT, Cambridge, MA 02139, USA}
\email{tabuada@math.mit.edu}
\urladdr{http://math.mit.edu/~tabuada}
\thanks{The author is very grateful to the organizers Nitu Kitchloo, Mona Merling, Jack Morava, Emily Riehl, and W. Stephen Wilson, for the kind invitation to present some of this work at the second Mid-Atlantic Topology Conference. The author was supported by~a~NSF~CAREER~Award.}

\date{\today}

\abstract{This survey covers some of the recent developments on noncommutative motives and their applications. Among other topics, we compute the additive invariants of relative cellular spaces and orbifolds; prove Kontsevich's semi-simplicity conjecture; prove a far-reaching noncommutative generalization of the Weil conjectures; prove Grothendieck's standard conjectures of type $C^+$ and $D$, Voevodsky's nilpotence conjecture, and Tate's conjecture, in several new cases; embed the (cohomological) Brauer group into secondary $K$-theory; construct a noncommutative motivic Gysin triangle; compute the  localizing $\bbA^1$-homotopy invariants of corner skew Laurent polynomial algebras and of noncommutative projective schemes; relate Kontsevich's category of noncommutative mixed motives to Morel-Voevodsky's stable $\bbA^1$-homotopy category, to Voevodsky's triangulated category of mixed motives, and to Levine's triangulated category of mixed motives; prove the Schur-finiteness conjecture for quadric fibrations over low-dimensional bases; and finally extend Grothendieck's theory of periods to the setting of dg categories.}
}

\maketitle

\smallskip

\vskip-\baselineskip

\vspace{-0.3cm}

\quad \quad \quad \quad \quad \quad \quad \quad \quad \quad \quad \quad \quad \quad \quad \quad \quad \quad \quad \,{\em {\small To Lily, for being by my side.}}
\section*{Introduction}
After the release of the monograph [Noncommutative Motives. With a preface by Yuri I. Manin. University Lecture Series {\bf 63}, American Mathematical Society, 2015], several important results on the theory of noncommutative motives have been established. The purpose of this survey, written for a broad mathematical audience, is to give a rigorous overview of some of these recent results. We will follow closely the notations, as well as the writing style, of the monograph \cite{book}. Therefore, we suggest the reader to have it at his/her desk while reading this survey. The monograph \cite{book} is divided into the following chapters:

\vspace{0.2cm}

\begin{itemize}
\item[] Chapter 1. Differential graded categories.
\item[] Chapter 2. Additive invariants.
\item[] Chapter 3. Background on pure motives.
\item[] Chapter 4. Noncommutative pure motives.
\item[] Chapter 5. Noncommutative (standard) conjectures.
\item[] Chapter 6. Noncommutative motivic Galois groups.
\item[] Chapter 7. Jacobians of noncommutative Chow motives.
\item[] Chapter 8. Localizing invariants.
\item[] Chapter 9. Noncommutative mixed motives.
\item[] Chapter 10. Noncommutative motivic Hopf dg algebras
\item[] Appendix A. Grothendieck derivators.
\end{itemize}

\vspace{0.2cm}

In this survey we cover some of the recent developments concerning the Chapters 2, 4, 5, 6, 8, and 9. These developments are described in Sections \ref{sec:additive}, \ref{sec:pure}, \ref{sec:conjectures}, \ref{sec:Galois}, \ref{sec:localizing}, and \ref{sec:NCmixed}, respectively. The final Section \ref{sec:periods}, entitled ``Noncommutative realizations and periods'', discusses a recent research subject which was not addressed in \cite{book}.
\subsection*{Preliminaries}
Throughout the survey, $k$ will denote a base field. We will assume the reader is familiar with the language of differential graded (=dg) categories; for a survey on dg categories, we invite the reader to consult Keller's ICM address \cite{ICM-Keller}. In particular, we will freely use the notions of {\em Morita equivalence} of dg categories (see \cite[\S1.6]{book}) and {\em smooth/proper dg category} in the sense of Kontsevich (see \cite[\S1.7]{book}). We will write $\dgcat(k)$ for the category of (small) dg categories and $\dgcat_\mathrm{sp}(k)$ for the full subcategory of smooth proper dg categories. Given a $k$-scheme $X$ (or more generally an algebraic stack $\cX$), we will denote by $\perf_\dg(X)$ the canonical dg enhancement of the category of perfect complexes $\perf(X)$; see \cite[Example~1.27]{book}.
\section{Additive invariants}\label{sec:additive}
Recall from \cite[\S2.3]{book} the construction of the {\em universal additive invariant} of dg categories $U\colon \dgcat(k) \to \Hmo_0(k)$. In \cite[\S2.4]{book} we described the behavior of $U$ with respect to semi-orthogonal decompositions, full exceptional collections, purely inseparable field extensions, central simple algebras, sheaves of Azumaya algebras, twisted flag varieties, nilpotent ideals, finite-dimensional algebras of finite global dimension, etc. In \S\ref{sub:relative}-\ref{sub:orbifolds} we describe the behavior of $U$ with respect to relative cellular spaces and orbifolds. As explained in \cite[Thm.~2.9]{book}, all the results in \S\ref{sub:relative}-\ref{sub:orbifolds} are {\em motivic} in the sense that they hold similarly for every additive invariant such as algebraic $K$-theory, mod-$n$ algebraic $K$-theory, Karoubi-Villamayor $K$-theory, nonconnective algebraic $K$-theory, homotopy $K$-theory, {\'e}tale $K$-theory, Hochschild homology, cyclic homology, negative cyclic homology, periodic cyclic homology, topological Hochschild homology, topological cyclic homology, topological periodic cyclic homology, etc. Consult \S\ref{sec:NCrealizations} for further examples of additive invariants.

\begin{notation} Given a $k$-scheme $X$ (or more generally an algebraic stack $\cX$), we will write $U(X)$ instead of $U(\perf_\dg(X))$.
\end{notation}
\subsection{Relative cellular spaces}\label{sub:relative}
A flat morphism of $k$-schemes $p\colon X \to Y$ is called an {\em affine fibration} of relative dimension $d$ if for every point $y\in Y$ there exists a Zariski open neighborhood $y \in V$ such that $X_V:=p^{-1}(V) \simeq Y \times \bbA^d$ with $p_V\colon X_V \to Y$ isomorphic to the projection onto the first factor. Following Karpenko \cite[Def.~6.1]{Karpenko}, a smooth projective $k$-scheme $X$ is called a {\em relative cellular space} if it admits a filtration by closed subschemes
\begin{equation*}
\varnothing = X_{-1} \hookrightarrow X_0 \hookrightarrow \cdots \hookrightarrow X_i \hookrightarrow \cdots \hookrightarrow X_{n-1} \hookrightarrow X_n=X
\end{equation*}
and affine fibrations $p_i\colon X_i \backslash X_{i-1} \to Y_i, 0 \leq i \leq n$, of relative dimension $d_i$ with $Y_i$ a smooth projective $k$-scheme. The smooth $k$-schemes $X_i \backslash X_{i-1}$ are called the {\em cells} and the smooth projective $k$-schemes $Y_i$ the {\em bases} of the cells.
\begin{example}[$\bbG_m$-schemes]\label{example:bb}
The celebrated  Bialynicki-Birula decomposition \cite{BB} provides a relative cellular space structure on smooth projective $k$-schemes equipped with a $\bbG_m$-action. In this case, the bases of the cells are given by the connected components of the fixed point locus.
This class of relative cellular spaces includes  the isotropic flag varieties considered by Karpenko in \cite{Karpenko} as well as the isotropic homogeneous spaces considered by Chernousov-Gille-Merkurjev~in~\cite{CGM}.
\end{example}
\begin{theorem}[{\cite[Thm.~2.7]{Gysin}}]\label{thm:cellular}
Given a relative cellular space $X$, we have an isomorphism $U(X) \simeq \bigoplus_{i=0}^n U(Y_i)$. 
\end{theorem}
Theorem \ref{thm:cellular} shows that the additive invariants of relative cellular spaces $X$ are completely determined by the basis $Y_i$ of the cells $X_i \backslash X_{i-1}$. Among other ingredients, its proof makes use of Theorem \ref{thm:Gysin}; consult \cite[\S9]{Gysin} for details.
\begin{example}[Kn\"orrer periodicity] 
Let $q=fg+q'$, where $f$, $g$, and $q'$, are
  forms of degrees $a>0$, $b>0$, and $a+b$, in disjoint sets
  of variables $(x_i)_{i=1,\ldots,m}$, $(y_j)_{j=1,\ldots,n}$, and
  $(z_l)_{l=1,\ldots,p}$, respectively. Such a decomposition holds, for example, in the case of isotropic quadratic forms $q$. Let us write $Q$ and $Q'$ for the projective
  hypersurfaces defined by $q$ and $q'$, respectively. Assume that $Q$ is smooth. Under this assumption, we have a $\bbG_m$-action on $Q$ given by $\lambda\cdot
  (\underline{x},\underline{y},\underline{z}):=(\lambda^b\underline{x},\lambda^{-a}\underline{y},\underline{z})$
  with fixed point locus $\bbP^{m-1}\amalg \bbP^{n-1}\amalg Q'$; this implies that $Q'$ is also smooth. By combining Theorem \ref{thm:cellular} and Example \ref{example:bb} with the fact that $U(\bbP^n)\simeq U(k)^{\oplus (n+1)}$ (see \cite[\S2.4.2]{book}), we obtain an induced isomorphism $U(Q)\simeq U(k)^{\oplus (m+n)}\oplus U(Q')$. Morally speaking, this shows that (modulo $k$) the additive invariants of $Q$ and $Q'$ are the same.
\end{example}
\subsection{Orbifolds}\label{sub:orbifolds}
Let $G$ be a finite group of order $n$ (we assume that $1/n \in k$), $\varphi$ the set of all cyclic subgroups of $G$, $\vphi$ a set of representatives of the conjugacy classes in $\varphi$, $X$ a smooth $k$-scheme equipped with a $G$-action, and $[X/G]$ the associated orbifold. As explained in \cite[\S3]{Orbifold}, the assignment $[V] \mapsto V\otimes_k -$, where $V$ stands for a $G$-representation, gives rise to an action of the representation ring $R(G)$ on $U([X/G])$. Given $\sigma \in \varphi$, let $e_\sigma$ be the unique idempotent of the $\bbZ[1/n]$-linearized representation ring $R(\sigma)_{1/n}$ whose image under all the restrictions $R(\sigma)_{1/n}~\to~R(\sigma')_{1/n}$, with $\sigma' \subsetneq \sigma$, is zero. 
The normalizer $N(\sigma)$ of $\sigma$ acts naturally on $[X^\sigma/\sigma]$ and hence on $U([X^\sigma/\sigma])$. By functoriality, this action restricts to the direct summand $e_\sigma U([X^\sigma/\sigma])_{1/n}$.
\begin{theorem}[{\cite[Thm.~1.1 and Cor.~1.6]{Orbifold}}]\label{thm:orbifold}
The following computations hold:
\begin{itemize}
\item[(i)] We have an induced isomorphism
\begin{equation}\label{eq:formula1}
U([X/G])_{1/n}\simeq \bigoplus_{\sigma \in \varphi\!/\!\sim} (e_\sigma U([X^\sigma/\sigma])_{1/n})^{N(\sigma)}
\end{equation}
in the $\bbZ[1/n]$-linearized (and idempotent completed) category $\Hmo_0(k)_{1/n}$.
\item[(ii)] If $k$ contains the $n^{\mathrm{th}}$ roots of unity, then \eqref{eq:formula1} reduces to an isomorphism
\begin{equation}\label{eq:formula2}
U([X/G])_{1/n} \simeq \bigoplus_{\sigma \in \varphi\!/\!\sim}(U(X^\sigma)_{1/n} \otimes_{\bbZ[1/n]} e_\sigma R(\sigma)_{1/n})^{N(\sigma)}\,,
\end{equation}
where $-\otimes_{\bbZ[1/n]}-$ stands for the canonical action of the category of finitely generated projective $\bbZ[1/n]$-modules on $\Hmo_0(k)_{1/n}$.
\item[(iii)] If $k$ contains the $n^{\mathrm{th}}$ roots of unity and $F$ is a field which contains the $n^{\mathrm{th}}$ roots of unity and $1/n \in F$, then we have induced isomorphisms
\begin{equation}\label{eq:formula3}
U([X/G])_F\simeq \bigoplus_{g \in G\!/\!\sim} U(X^g)_F^{C(g)} \simeq (\bigoplus_{g \in G} U(X^g)_F)^G
\end{equation}
in the category $\Hmo_0(k)_F$, where $C(g)$ stands for the centralizer of $g$.
\end{itemize}
Moreover, \eqref{eq:formula2}-\eqref{eq:formula3} are isomorphisms of (commutative) monoids.
\end{theorem}
Roughly speaking, Theorem \ref{thm:orbifold} shows that the additive invariants of orbifolds can be computed using solely ``ordinary'' schemes. 
\begin{example}[McKay correspondence]
In many cases, the dg category $\perf_\dg([X/G])$ is known to be Morita equivalent 
to $\perf_\dg(Y)$ for a crepant 
  resolution $Y$ of the (singular) geometric quotient $X/\!\!/G$; see \cite{BezKal2,BKR,KV,
    Kawamata}. This is generally referred to as the ``McKay correspondence''.
Whenever it holds, we can replace $[X/G]$ by $Y$ in the formulas \eqref{eq:formula1}-\eqref{eq:formula3}. Here is an illustrative example (with $k$ algebraically closed): the cyclic group $G=C_2$ acts on any abelian surface $S$ by the
  involution $a \mapsto -a$ and the Kummer surface $\mathrm{Km}(S)$ is defined as
the blow-up of $S/\!\!/C_2$ at its $16$ singular points.
In this case, the dg category $\perf_\dg([S/C_2])$ is Morita equivalent to 
  $\perf_\dg(\mathrm{Km}(S))$. Consequently, Theorem \ref{thm:orbifold}(ii) leads to an isomorphism:
\begin{equation}\label{eq:McKay}
U(\operatorname{Km}(S))_{1/2}\simeq U(S)_{1/2}^{C_2} \oplus U(k)_{1/2}^{\oplus 16}\,.
\end{equation}
Note that since the Kummer surface is Calabi-Yau, the category $\perf(\operatorname{Km}(S))$ does {\em not} admit non-trivial semi-orthogonal decompositions. This shows that the isomorphism \eqref{eq:McKay} is {\em not} induced from a semi-orthogonal decomposition.
\end{example}
\begin{corollary}[Algebraic $K$-theory]
If $k$ contains the $n^{\mathrm{th}}$ roots of unity, then we have the following isomorphism of $\bbZ$-graded commutative $\bbZ[1/n]$-algebras:
\begin{equation}\label{eq:formula4}
 K_\ast([X/G])_{1/n} \simeq \bigoplus_{\sigma \in \varphi\!/\!\sim}(K_\ast(X^\sigma)_{1/n} \otimes_{\bbZ[1/n]} e_\sigma R(\sigma)_{1/n})^{N(\sigma)}\,.
\end{equation}
\end{corollary}
The formula \eqref{eq:formula4} was originally established by Vistoli in \cite[Thm.~1]{Vistoli}. Among other ingredients, Vistoli's proof makes essential use of d\'evissage. The proof of Theorem \ref{thm:orbifold}, and hence of \eqref{eq:formula4}, is not only different but moreover avoids the use of d\'evissage; consult \cite[\S6]{Orbifold} for details.
\begin{corollary}[Cyclic homology]
If $k$ contains the $n^{\mathrm{th}}$ roots of unity, then we have the following isomorphisms of $\bbZ$-graded commutative $k$-algebras: 
\begin{equation}\label{eq:formula5}
HC_\ast([X/G]) \simeq \bigoplus_{g\in G\!/\!\sim} HC_\ast(X^g)^{C(g)} \simeq (\bigoplus_{g\in G} HC_\ast(X^g))^G\,.
\end{equation}
\end{corollary}
The formula \eqref{eq:formula5} was originally established by Baranovsky in \cite[Thm.~1.1]{Baranovsky}. Baranovsky's proof is very specific to cyclic homology. In constrast, the proof of Theorem \ref{thm:orbifold}, and hence of \eqref{eq:formula5}, avoids all the specificities of cyclic homology and is moreover quite conceptual; consult \cite[\S6]{Orbifold} for details.
\begin{corollary}[Topological periodic cyclic homology]
Let $k$ be a perfect field of characteristic $p>0$, $W(k)$ the associated ring of $p$-typical Witt vectors, and $K:=W(k)[1/p]$ the fraction field of $W(k)$. If $k$ contains the $n^{\mathrm{th}}$ roots of unity, then we have the following isomorphisms of $\bbZ/2$-graded commutative $K$-algebras: 
\begin{equation}\label{eq:formula6}
TP_\ast([X/G])_{1/p} \simeq \bigoplus_{g\in G\!/\!\sim} TP_\ast(X^g)_{1/p}^{C(g)} \simeq (\bigoplus_{g\in G} TP_\ast(X^g)_{1/p})^G\,.
\end{equation}
\end{corollary}
To the best of the author's knowledge, the formula \eqref{eq:formula6} is new in the literature; consult \cite[\S1]{Orbifold} for further corollaries of Theorem \ref{thm:orbifold}.
\subsubsection{Twisted analogues}\label{sec:twisted}
Given a sheaf of Azumaya algebras $\cF$ over $[X/G]$, \ie a $G$-equivariant sheaf of Azumaya algebras over $X$, all the computations of Theorem \ref{thm:orbifold} admit $\cF$-twisted analogues; consult \cite[Thm.~1.27 and Cor.~1.29]{Orbifold} for details.
\section{Noncommutative pure motives}\label{sec:pure}
In \S\ref{sub:recollections} we recall the definition of the different categories of noncommutative pure motives. Subsections \S\ref{sub:Brauer}-\ref{sub:rigidity} are devoted to three structural properties of these categories (relation with the Brauer group, semi-simplicity and rigidity). In \ref{sub:zeta} we prove to a far-reaching noncommutative generalization of the Weil conjectures; see Theorem \ref{thm:zeta}. Finally, in \S\ref{sub:equivariant} we describe some of the equivariant analogues of the theory of noncommutative pure motives.
\subsection{Recollections}\label{sub:recollections}
Recall from \cite[\S4.1]{book} that the category of {\em noncommutative Chow motives $\NChow(k)$} is defined as the idempotent completion of the full subcategory of $\Hmo_0(k)$ consisting of the objects $U(\cA)$, with $\cA$ a smooth proper dg category. By construction, this category is additive, rigid symmetric monoidal, and comes equipped with a symmetric monoidal functor $U\colon \dgcat_{\mathrm{sp}}(k) \to \NChow(k)$. Moreover, given smooth proper dg categories $\cA$ and $\cB$, we have isomorphisms:
\begin{equation}\label{eq:star}
\Hom_{\NChow(k)}(U(\cA),U(\cB)) \simeq K_0(\cD_c(\cA^\op \otimes \cB))=: K_0(\cA^\op \otimes \cB)\,.
\end{equation}
Given a rigid symmetric monoidal category $(\cC, \otimes, {\bf 1})$, consider the $\otimes$-ideal
$$\otimes_{\mathrm{nil}}(a,b) := \{ f \in \Hom_\cC(a,b) \,|\,f^{\otimes n}=0\,\,\mathrm{for}\,\,\mathrm{some}\,\,n\gg 0\}\,.$$
Recall from \cite[\S4.4]{book} that the category of {\em noncommutative $\otimes$-nilpotent motives $\NVoev(k)$} is defined as the idempotent completion of the quotient $\NChow(k)/\otimes_{\mathrm{nil}}$.

As explained in \cite[\S4.5]{book}, periodic cyclic homology gives rise to an additive symmetric monoidal functor $HP^\pm\colon \NChow(k) \to \mathrm{Vect}_{\bbZ/2}(k)$, with values in the category of finite-dimensional $\bbZ/2$-graded $k$-vector spaces. Recall from {\em loc. cit.} that the category of {\em noncommutative homological motives $\NHom(k)$} is defined as the idempotent completion of the quotient $\NChow(k)/\mathrm{Ker}(HP^\pm)$.

Given a rigid symmetric monoidal category $(\cC, \otimes, {\bf 1})$, consider the $\otimes$-ideal
$$
\cN(a,b) := \{ f \in \Hom_\cC(a,b)\,|\,\forall g \in \Hom_\cC(b,a)\,\,\mathrm{we}\,\,\mathrm{have}\,\,\mathrm{tr}(g\circ f)=0\}\,,
$$
where $\mathrm{tr}(g\circ f)$ stands for the categorical trace of the endomorphism $g\circ f$. Recall from \cite[\S4.6]{book} that the category of {\em noncommutative numerical motives $\NNum(k)$} is defined as the idempotent completion of the quotient $\NChow(k)/\cN$.
\subsection{Relation with the Brauer group}\label{sub:Brauer}
Let $\mathrm{Br}(k)$ be the Brauer group of the base field $k$. Given a central simple $k$-algebra $A$, we write $[A]$ for its Brauer class.
\begin{example}[Local fields]
A local field $k$ is isomorphic to $\bbR$, to $\bbC$, to a finite field extension of $\bbQ_p$, or to a finite field extension of $\bbF_p(\!(t)\!)$. Thanks to local class field theory, we have $\mathrm{Br}(\bbR)\simeq \bbZ/2$, $\mathrm{Br}(\bbC)=0$, and $\mathrm{Br}(k)\simeq \bbQ/\bbZ$ in all the remaining cases. Moreover, every element of $\mathrm{Br}(k)$ can be represented by a cyclic $k$-algebra.
\end{example}
Recall from \cite[\S2.4.4]{book} that we have the following equivalence
\begin{equation}\label{eq:equivalence-1}
[A]=[B] \Leftrightarrow U(A) \simeq U(B)
\end{equation}
for any two central simple $k$-algebras $A$ and $B$. Intuitively speaking, \eqref{eq:equivalence-1} shows that the Brauer class $[A]$ and the noncommutative Chow motive $U(A)$ contain exactly the same information. Let $K_0(\NChow(k))$ be the Grothendieck ring of the additive symmetric monoidal category $\NChow(k)$. Given a central simple $k$-algebra $A$, we write $[U(A)]$ for the Grothendieck class of $U(A)$. The (proof of the) next result is contained in \cite[Thm.~6.12]{Secondary}\cite[Thm.~1.3]{SecondaryII}:
\begin{theorem}\label{thm:new}
Given central simple $k$-algebras $A$ and $B$, we have the equivalence:
\begin{equation}\label{eq:equivalence-2}
U(A)\simeq U(B) \Leftrightarrow [U(A)] \simeq [U(B)]\,.
\end{equation}
\end{theorem}
Roughly speaking, Theorem \ref{thm:new} shows that the noncommutative Chow motives of central simple $k$-algebras are insensitive to the Grothendieck group relations. By combining the equivalences \eqref{eq:equivalence-1} and \eqref{eq:equivalence-2}, we obtain the following result:
\begin{corollary}\label{cor:new}
The following map is injective:
\begin{eqnarray*}
\mathrm{Br}(k) \too K_0(\NChow(k)) && [A] \mapsto [U(A)]\,.
\end{eqnarray*}
\end{corollary}
Consult \S\ref{sub:secondary} for some applications of Corollary \ref{cor:new} to secondary $K$-theory.
\begin{remark}[Generalizations]\label{rk:generalization1}
Theorem \ref{thm:new} and Corollary \ref{cor:new} hold more generally with $k$ replaced by a base $k$-scheme $X$. Furthermore, instead of the Brauer group $\mathrm{Br}(X)$, we can consider the second \'etale cohomology group\footnote{As proved by Gabber \cite{Gabber} and de Jong \cite{deJong}, in the case where $X$ admits an ample line bundle (\eg\ $X$ affine), the Brauer group $\mathrm{Br}(X)$ may be identified with the torsion subgroup~of~$H^2_{\mathrm{et}}(X,\bbG_m)$.} $H^2_{\mathrm{et}}(X,\bbG_m)$; consult \cite{Secondary,SecondaryII} for details. In the case of an affine cone over a smooth irreducible plane complex curve of degree $\geq 4$, the latter \'etale cohomology group contains non-torsion classes. The same phenomenon occurs, for example, in the case of Mumford's (celebrated) singular surface \cite[page~75]{Mumford}; see \cite[Example~1.32]{SecondaryII}.
\end{remark}
\begin{remark}[Jacques Tits' motivic measure]
The Grothendieck ring of varieties $K_0\mathrm{Var}(k)$, introduced in a letter from Grothendieck to Serre in the sixties, is defined as the quotient of the free abelian group on the set of isomorphism classes of $k$-schemes by the ``cut-and-paste'' relations. Although very important, the structure of this ring still remains poorly understood. Among other ingredients, Theorem \ref{thm:new} was used in the construction of a new motivic measure $\mu_T$ entitled {\em Tits motivic measure}; consult \cite{Tits} for details. This new motivic measure led to the proof of several new structural properties of $K_0\mathrm{Var}(k)$. For example, making use of $\mu_T$, it was proved in {\em loc. cit.} that two quadric hypersurfaces (or more generally involution varieties), associated to quadratic forms of degree $6$, have the same Grothendieck class if and only if they are isomorphic. In the same vein, it was proved in {\em loc. cit.} that two products of conics have the same Grothendieck class if and only if they are isomorphic; this refines a previous result of Koll\'ar \cite{Kollar}.
\end{remark}
\subsection{Semi-simplicity}\label{sub:semi-simple}
Let $F$ be a field of characteristic zero. The following result is obtained by combining \cite[Thm.~4.27]{book} with \cite[Thm.~1.1]{positive}:
\begin{theorem}\label{thm:semi-simple}
The category $\NNum(k)_F$ is abelian semi-simple.
\end{theorem}
Assuming certain (polarization) conjectures, Kontsevich conjectured in his seminal talk \cite{IAS} that the category $\NNum(k)_F$ was abelian semi-simple. Theorem \ref{thm:semi-simple} not only proves this conjecture but moreover shows that Kontsevich's insight holds unconditionally. Let $\Num(k)_F$ be the (classical) category of numerical motives; see \cite[\S4]{Andre}. The next result is obtained by combining \cite[Rk.~4.32]{book} with \cite[Cor.~1.2]{positive}:
\begin{corollary}\label{cor:semi-simple}
The category $\Num(k)_F$ is abelian semi-simple.
\end{corollary}
Assuming certain (standard) conjectures, Grothendieck conjectured in the sixties that the category $\Num(k)_F$ was abelian semi-simple. This conjecture was proved unconditionally by Jannsen \cite{Jannsen} in the nineties using \'etale cohomology. Corollary \ref{cor:semi-simple} provides us with an alternative proof of Grothendieck's conjecture.
\subsubsection{Numerical Grothendieck group}
The Grothendieck group $K_0(\cA)$ of a proper dg category $\cA$ comes equipped with the following Euler bilinear pairing:
\begin{eqnarray*}\label{eq:pairing}
\chi \colon  K_0(\cA) \times K_0(\cA) \too \bbZ && ([M],[N]) \mapsto \sum_n (-1)^n \mathrm{dim}_k \Hom_{\cD_c(\cA)}(M,N[n])\,.
\end{eqnarray*}
This bilinear pairing is, in general, not symmetric neither skew-symmetric. Nevertheless, when $\cA$ is moreover smooth the associated left and right kernels of $\chi$ agree; see \cite[Prop.~4.24]{book}. Consequently, under these assumptions on $\cA$, we have a well-defined {\em numerical Grothendieck group} $K_0(\cA)/_{\!\!\sim \mathrm{num}}:=K_0(\cA)/\mathrm{Ker}(\chi)$. Following \cite[Thm.~4.26]{book}, given smooth proper dg categories $\cA$ and $\cB$, we have isomorphisms:
\begin{equation}\label{eq:natural-iso}
\Hom_{\NNum(k)}(U(\cA),U(\cB)) \simeq K_0(\cA^\op \otimes \cB)/\mathrm{Ker}(\chi)\,.
\end{equation}
The next result, whose proof makes use of Theorem \ref{thm:semi-simple}, is obtained by combining \cite[Thm.~1.2]{Separable} with \cite[Thm.~6.2]{positive}:
\begin{theorem}\label{thm:free}
$K_0(\cA)/_{\!\!\sim \mathrm{num}}$ is a finitely generated free abelian group.
\end{theorem}
Given a smooth proper $k$-scheme $X$, let us write $\cZ^\ast(X)/_{\!\!\sim \mathrm{num}}$ for the (graded) group of algebraic cycles on $X$ up to numerical equivalence. By combining Theorem \ref{thm:free} with the Hirzebruch-Riemann-Roch theorem, we obtain the following result:
\begin{corollary}
$\cZ^\ast(X)/_{\!\!\sim \mathrm{num}}$ is a finitely generated free abelian (graded) group.
\end{corollary}
\subsection{Rigidity}\label{sub:rigidity}
Recall that a field extension $l/k$ is called {\em primary} if the algebraic closure of $k$ in $l$ is purely inseparable over $k$. When $k$ is algebraically closed, every field extension $l/k$ is primary.
\begin{theorem}[{\cite[Thm.~2.1(i)]{rigidity}}]\label{thm:rigidity}
Given a primary field extension $l/k$ and a field $F$ of characteristic zero, the base-change functor $-\otimes_k l: \NNum(k)_F \to \NNum(l)_F$ is fully-faithful. The same holds integrally when $k$ is algebraically closed. 
\end{theorem}
Intuitively speaking, Theorem \ref{thm:rigidity} shows that the theory of noncommutative numerical motives is ``rigid'' under base-change along primary field extensions. Alternatively, thanks to the isomorphisms \eqref{eq:natural-iso}, Theorem \ref{thm:rigidity} shows that the numerical Grothendieck group is ``rigid'' under primary field extensions. The commutative counterpart, resp. mixed analogue, of Theorem \ref{thm:rigidity} was established by Kahn in \cite[Prop.~5.5]{Kahn2}, resp. is provided by Theorem \ref{thm:rigidity2}.
\begin{remark}[Extra functoriality]
Let $l/k$ be a primary field extension. As proved in \cite[Thm.~2.3]{rigidity}, Theorems \ref{thm:semi-simple} and \ref{thm:rigidity} imply that the base-change functor admits a left=right adjoint. Without the assumption that the field extension $l/k$ is primary, such an adjoint functor does {\em not} exists in general; consult \cite[Rk.~2.4]{rigidity} for details.
\end{remark}
\subsection{Zeta functions of endomorphisms}\label{sub:zeta}
Let $N\!\!M \in \NChow(k)_\bbQ$ be a noncommutative Chow motive and $f$ an endomorphism of $N\!\!M$. Following Kahn \cite[Def.~3.1]{Zeta}, the {\em zeta function of $f$} is defined as the following formal power series
\begin{equation}\label{eq:zeta}
Z(f;t):= \mathrm{exp}\left(\sum_{n \geq 1} \mathrm{tr}(f^{\circ n}) \frac{t^n}{n}\right) \in \bbQ\llbracket t \rrbracket\,,
\end{equation}
where $f^{\circ n}$ stands for the composition of $f$ with itself $n$-times, $\mathrm{tr}(f^{\circ n}) \in \bbQ$ stands for the categorical trace of $f^{\circ n}$, and $\mathrm{exp}(t):=\sum_{m \geq 0} \frac{t^m}{m!} \in \bbQ\llbracket t\rrbracket$. 
\begin{remark}
When $N\!\!M = U(\cA)_\bbQ$ and $f=[\mathrm{B}]_\bbQ$, with $\mathrm{B} \in \cD_c(\cA^\op \otimes \cA)$ a dg $\cA\text{-}\cA$-bimodule (see \S\ref{sub:recollections}), we have the following computation
\begin{equation}\label{eq:integers}
\mathrm{tr}(f^{\circ n}) =[HH(\cA; \underbrace{\mathrm{B}\otimes^{\bf L}_\cA \cdots \otimes^{\bf L}_\cA \mathrm{B}}_{n\text{-}\text{times}})] \in K_0(k) \simeq \bbZ\,,
\end{equation}
where $HH(\cA; \mathrm{B}\otimes^{\bf L}_\cA \cdots \otimes^{\bf L}_\cA \mathrm{B})$ stands for the Hochschild homology of $\cA$ with coefficients in $\mathrm{B}\otimes^{\bf L}_\cA \cdots \otimes^{\bf L}_\cA \mathrm{B}$; see \cite[Prop.~2.26]{book}. Intuitively speaking, the integer \eqref{eq:integers} is the ``number of fixed points'' of the dg $\cA\text{-}\cA$-bimodule $\mathrm{B}\otimes^{\bf L}_\cA \cdots \otimes^{\bf L}_\cA \mathrm{B}$.
\end{remark}
\begin{example}[Zeta function]\label{ex:schemes}
Let $k=\bbF_q$ be a finite field, $X$ a smooth proper $k$-scheme, and $\mathrm{Fr}$ the geometric Frobenius. When $\cA= \perf_\dg(X)$ and $\mathrm{B}$ is the dg bimodule associated to the pull-back dg functor $\mathrm{Fr}^\ast\colon \perf_\dg(X) \to \perf_\dg(X)$, \eqref{eq:integers} reduces to $[HH(X;\Gamma_{\mathrm{Fr}^{\circ n}})]=\langle \Delta \cdot \Gamma_{\mathrm{Fr}^{\circ n}}\rangle= |X(\bbF_{q^n})|$. Consequently, \eqref{eq:zeta} reduces to the (classical) zeta function $Z_X(t):= \mathrm{exp}(\sum_{n \geq 1} |X(\bbF_{q^n})| \frac{t^n}{n})$ of $X$.
\end{example}
\begin{remark}[Witt vectors]
Recall from \cite{Hazewinkel} the definition of the ring of (big) Witt vectors $\mathrm{W}(\bbQ)=(1 + t \bbQ\llbracket t \rrbracket, \times, \ast)$. Since the leading term of \eqref{eq:zeta} is equal to $1$, the zeta function $Z(f;t)$ of $f$ belongs to $\mathrm{W}(\bbQ)$. Moreover, given endomorphisms $f$ and $f'$ of noncommutative Chow motives $N\!\!M$ and $N\!\!M'$, we have $Z(f\oplus f';t)=Z(f;t) \times Z(f';t)$ and $Z(f\otimes f'; t) = Z(f;t) \ast Z(f';t)$ in $\mathrm{W}(\bbQ)$.
\end{remark}
Let $B = \prod_i B_i$ be a finite-dimensional semi-simple $\bbQ$-algebra, $Z_i$ the center of $B_i$, $\delta_i$ for the degree $[Z_i:\bbQ]$, and $d_i$ the index $[B_i: Z_i]^{1/2}$. Given a unit $b \in B^\times$, its {\em $i^{\mathrm{th}}$ reduced norm} $\mathrm{Nrd}_i(b) \in \bbQ$ is defined as the composition $(\mathrm{N}_{Z_i/\bbQ}\circ \mathrm{Nrd}_{B_i/Z_i})(b_i)$. 

Let $N\!\!M \in \NChow(k)_\bbQ$ be a noncommutative Chow motive. Thanks to Theorem \ref{thm:semi-simple}, $B:=\mathrm{End}_{\NNum(k)_\bbQ}(N\!\!M)$ is a finite-dimensional semi-simple $\bbQ$-algebra; let us write $e_i \in B$ for the central idempotent corresponding to the summand $B_i$. Given an invertible endomorphism $f$ of  $N\!\!M$, its {\em determinant $\mathrm{det}(f) \in \bbQ$} is defined as the following product $\prod_i \mathrm{Nrd}_i(f)^{\mu_i}$, where $\mu_i :=\frac{\mathrm{tr}(e_i)}{\delta_i d_i}$.
\begin{theorem}[{\cite[Thm.~5.8]{positive}}]\label{thm:zeta}
\begin{itemize}
\item[(i)] The series $Z(f;t) \in \bbQ \llbracket t \rrbracket$ is {\em rational}, \ie $Z(f;t)=\frac{p(t)}{q(t)}$ with $p(t), q(t) \in \bbQ[t]$. Moreover, $\mathrm{deg}(q(t)) - \mathrm{deg}(p(t))= \mathrm{tr}(\id_{N\!\!M})$.  
\item[(ii)] When $f$ is invertible, we have the following functional equation:
$$ Z(f^{-1};t^{-1}) = (-t)^{\mathrm{tr}(\id_{N\!\!M})} \mathrm{det}(f) Z(f;t)\,.$$
\end{itemize}
\end{theorem}
\begin{corollary}[Weil conjectures]\label{cor:zeta}
Let $k=\bbF_q$ be a finite field, $X$ a smooth proper $k$-scheme $X$ of dimension $d$, and $\mathrm{E}:=\langle \Delta \cdot \Delta \rangle \in \bbZ$ the self-intersection number of the diagonal $\Delta$ of $X \times X$. 
\begin{itemize}
\item[(i)] The zeta function $Z_X(t)$ of $X$ is rational. Moreover, $\mathrm{deg}(q(t))-\mathrm{deg}(p(t))= \mathrm{E}$.
\item[(ii)] We have the following functional equation $Z_X(\frac{1}{q^d t}) = \pm t^{\mathrm{E}} q^{\frac{d}{2} \mathrm{E}} Z_X(t)$.
\end{itemize}
\end{corollary}
Weil conjectured\footnote{Weil conjectured also that the zeta function $Z_X(t)$ of $X$ satisfied an analogue of the Riemann hypothesis. This conjecture was proved by Deligne \cite{Deligne-IHES} using, among other tools, Lefschetz pencils.} in \cite{Weil} that the zeta function $Z_X(t)$ of $X$ was rational and that it satisfied a functional equation. These conjectures were proved independently by Dwork \cite{Dwork} and Grothendieck \cite{Grothendieck-Bourbaki} using $p$-adic analysis and \'etale cohomology, respectively. Corollary \ref{cor:zeta} provides us with an alternative proof of the Weil conjectures; see \cite[Cor.~5.12]{positive}. Moreover, Theorem \ref{thm:zeta} proves a far-reaching noncommutative generalization of the Weil conjectures. 
\subsection{Equivariant noncommutative motives}\label{sub:equivariant}
Let $G$ be a finite group of order $n$ (we assume that $1/n \in k$). Recall from \cite[Def.~4.1]{Equivariant} the definition of a {\em $G$-action} on a dg category $\cA$. Given a $G$-action $G \circlearrowright \cA$, we have an associated dg category $\cA^G$ of $G$-equivariant objects. From a topological viewpoint, $\cA^{G}$ may be understood as the ``homotopy fixed points'' of the $G$-action on $\cA$. Here are two examples:
\begin{example}[$G$-schemes]
Given a $G$-scheme $X$, the dg category $\perf_\dg(X)$ inherits a $G$-action. In this case, the dg category $\perf_\dg(X)^G$ is Morita equivalent to the dg category of $G$-equivariant perfect complexes $\perf^{G}_\dg(X)=\perf_\dg([X/G])$.
\end{example}
\begin{example}[Cohomology classes]
Given a cohomology class $[\alpha] \in H^2(G,k^\times)$, the dg category $k$ inherits a $G$-action $G \circlearrowright_\alpha k$. In this case, the dg category of $G$-equivariant objects is Morita equivalent to the twisted group algebra $k_\alpha[G]$. 
\end{example}
Let $\dgcat^{G}(k)$ be the category of (small) dg categories equipped with a $G$-action, and $\dgcat^{G}_{\mathrm{sp}}(k)$ the full subcategory of smooth proper dg categories. As explained in \cite[\S5]{Equivariant}, the category $\NChow(k)$ admits a $G$-equivariant counterpart $\NChow^{G}(k)$. Recall from {\em loc. cit.} that the latter category is additive, rigid symmetric monoidal, and comes equipped with a symmetric monoidal functor $U^{G}\colon \dgcat_{\mathrm{sp}}^{G}(k) \to \NChow^{G}(k)$. Moreover, we have isomorphisms
$$\mathrm{Hom}_{\NChow^{G}(k)}(U^{G}(G \circlearrowright \cA), U^{G}(G \circlearrowright \cB)) \simeq K_0^{G}(\cA^\op \otimes \cB)\,,$$ 
where the right-hand side stands for the $G$-equivariant Grothendieck group. In particular, the ring of endomorphisms of the $\otimes$-unit $U^{G}(G \circlearrowright_1 k)$ agrees with the representation ring\footnote{Recall that when $k=\bbC$ and $G$ is abelian, we have an isomorphism $R(G)\simeq \bbZ[\widehat{G}]$.} $R(G)$. Let us write $I$ for the augmentation ideal associated to the rank homomorphism $R(G) \twoheadrightarrow \bbZ$.
\subsubsection{Relation with equivariant Chow motives}
Making use of Edidin-Graham's work \cite{EG} on equivariant intersection theory, Laterveer \cite{Laterveer}, and Iyer and M\"uller-Stach \cite{Iyer-Muller}, extended the theory of Chow motives to the $G$-equivariant setting. In particular, they constructed a category of $G$-equivariant Chow motives $\Chow^{G}(k)$ and a (contravariant) symmetric monoidal functor $\mathfrak{h}^{G}\colon \mathrm{SmProj}^{G}(k) \to \Chow^{G}(k)$, defined on smooth projective $G$-schemes.
\begin{theorem}[{\cite[Thm.~8.4]{Equivariant}}]\label{thm:bridge-equivariant}
There exists a $\bbQ$-linear, fully-faithful, symmetric monoidal $\Phi_\bbQ^{G}$ making the following diagram commute
\begin{equation*}
\xymatrix{
\mathrm{SmProj}^{G}(k) \ar[rrr]^-{X \mapsto G \circlearrowright \perf_\dg(X)} \ar[d]_-{\mathfrak{h}^{G}(-)_\bbQ} &&& \dgcat_{\mathrm{sp}}^{G}(k) \ar[d]^-{U^{G}(-)_\bbQ} \\
\Chow^{G}(k)_\bbQ \ar[d] &&& \NChow^{G}(k)_\bbQ \ar[d]^-{(-)_{I_\bbQ}} \\
\Chow^{G}(k)_\bbQ/_{\!-\otimes \bbQ(1)} \ar[rrr]_-{\Phi_\bbQ^{G}} &&& \NChow^{G}(k)_{\bbQ, I_\bbQ}\,,
}
\end{equation*}
where $\Chow^{G}(k)_\bbQ/_{\!-\otimes \bbQ(1)}$ stands for the orbit category of $\Chow^{G}(k)_\bbQ$ with respect to the $G$-equivariant Tate motive $\bbQ(1)$ (see \cite[\S4.2]{book}), and $(-)_{I_\bbQ}$ for the localization functor associated to the augmentation ideal $I_\bbQ$.
\end{theorem}
Roughly speaking, Theorem \ref{thm:bridge-equivariant} shows that in order to compare the equivariant commutative world with the equivariant noncommutative world, we need to ``$\otimes$-trivialize'' the $G$-equivariant Tate motive $\bbQ(1)$ on one side and to localize at the augmentation ideal $I_\bbQ$ on the other side. Only after these two reductions, the equivariant commutative world embeds fully-faithfully into the equivariant noncommutative world. As illustrated in \S\ref{sub:exceptional}, this shows that the $G$-equivariant Chow motive $\mathfrak{h}^{G}(X)_\bbQ$ and the $G$-equivariant noncommutative Chow motive $U^{G}(G \circlearrowright \perf_\dg(X))$ contain (important) independent information~about~$X$.
\subsubsection{Full exceptional collections}\label{sub:exceptional}
Let $X$ be a smooth projective $G$-scheme. In order to study it, we can proceed into two distinct directions. On one direction, we can associate to $X$ its $G$-equivariant Chow motive $\mathfrak{h}^{G}(X)_\bbQ$. On another direction, we can associate to $X$ the $G$-action $G \circlearrowright\perf_\dg(X)$. The following result, whose proof makes use of Theorem \ref{thm:bridge-equivariant}, relates these two distinct directions of study:
\begin{theorem}[{\cite[Thm.~1.2]{Equivariant}}]\label{thm:exceptional}
If the category $\perf(X)$ admits a full exceptional collection $(\cE_1, \ldots, \cE_n)$ of $G$-invariant objects ($\neq$ $G$-equivariant objects), then there exists a choice of integers $r_1, \ldots, r_n \in  \{0, \ldots, \mathrm{dim}(X)\}$ such that 
\begin{equation}\label{eq:decomp-motivic}
\mathfrak{h}^{G}(X)_\bbQ \simeq \bbL^{\otimes r_1} \oplus \cdots \oplus \bbL^{\otimes r_n}\,,
\end{equation}
where $\bbL\in \Chow^{G}(k)_\bbQ$ stands for the $G$-equivariant Lefschetz motive.
\end{theorem}
Theorem \ref{thm:exceptional} can be applied, for example, to any $G$-action on projective spaces, quadrics, Grassmannians, etc; consult \cite[Examples 9.9-9.11]{Equivariant} for details. Morally speaking, Theorem \ref{thm:exceptional} shows that the existence of a full exceptional collection of $G$-invariant objects completely determines the $G$-equivariant Chow motive $\mathfrak{h}^{G}(X)_\bbQ$. In particular, $\mathfrak{h}^{G}(X)_\bbQ$ loses all the information about the $G$-action on $X$. In contrast, as explained in \cite[Rmk.~9.4 and Prop.~9.8]{Equivariant}, the $G$-invariant objects $\cE_1, \ldots, \cE_n$ yield (non-trivial) cohomology classes $[\alpha_1], \ldots, [\alpha_n] \in H^2(G, k^\times)$ such that 
\begin{equation}\label{eq:motivic-decomp-NC}
U^{G}(G \circlearrowright\perf_\dg(X))\simeq U^{G}(G \circlearrowright_{\alpha_1} k) \oplus \cdots \oplus U^{G}(G \circlearrowright_{\alpha_n} k)\,.
\end{equation}
Taking into account \eqref{eq:decomp-motivic}-\eqref{eq:motivic-decomp-NC}, the $G$-equivariant Chow motive $\mathfrak{h}^{G}(X)_\bbQ$ and the $G$-equivariant noncommutative Chow motive $U^{G}(G \circlearrowright \perf_\dg(X))$ should be considered as complementary. While the former keeps track of the Tate twists but not of the $G$-action, the latter keeps track of the $G$-action but not of the Tate twists.
\section{Noncommutative (standard) conjectures}\label{sec:conjectures}
In \S\ref{sub:recollections2} we recall some important conjectures of Grothendieck, Voevodsky, and Tate. Subsection \S\ref{sub:counterparts} is devoted to their noncommutative counterparts. As a first application of the noncommutative viewpoint, we prove that the original conjectures of Grothendieck, Voevodsky, and Tate, are invariant under homological projective duality. This leads to a proof of these original conjectures in several new cases. As a second application, we extend the original conjectures from schemes to algebraic stacks and prove them in the case of ``low-dimensional'' orbifolds.
\subsection{Recollections}\label{sub:recollections2}
Let $k$ be a base field of characteristic zero. Given a smooth proper $k$-scheme $X$ and a Weil cohomology theory $H^\ast$, let us write $\pi_X^n$ for the $n^{\mathrm{th}}$ K\"unneth projector of $H^\ast(X)$, $Z^\ast(X)_\bbQ$ for the $\bbQ$-vector space of algebraic cycles on $X$, and $Z^\ast(X)_\bbQ/_{\!\sim \mathrm{nil}}$,  $Z^\ast(X)_\bbQ/_{\!\sim \mathrm{hom}}$, and $Z^\ast(X)_\bbQ/_{\!\sim \mathrm{num}}$, for the quotient of $Z^\ast(X)_\bbQ$ with respect to the smash-nilpotence, homological, and numerical equivalence relation, respectively. Recall from \cite[\S3.0.8-3.0.11]{book} that:
\begin{itemize}
\item[(i)] The Grothendieck's standard conjecture\footnote{The standard conjecture of type $C^+$ is also known as the {\em sign conjecture}. If the even K\"unneth projector $\pi_X^+$ is algebraic, then the odd K\"unneth projector $\pi^-_X:=\sum_n \pi_X^{2n+1}$ is also algebraic.} of type $C^+$, denoted by $C^+(X)$, asserts that the even K\"unneth projector $\pi^+_X:=\sum_n \pi^{2n}_X$ is algebraic.
\item[(ii)] The Grothendieck's standard conjecture of type $D$, denoted by $D(X)$, asserts that $Z^\ast(X)_\bbQ/_{\!\sim \mathrm{hom}}=Z^\ast(X)_\bbQ/_{\!\sim \mathrm{num}}$.
\item[(iii)] The Voevodsky's nilpotence conjecture $V(X)$ (which implies Grothendieck's conjecture $D(X)$) asserts that $Z^\ast(X)_\bbQ/_{\!\sim \mathrm{nil}}=Z^\ast(X)_\bbQ/_{\!\sim \mathrm{num}}$.
\item[(iv)] The Schur-finiteness conjecture\footnote{Consult \S\ref{sec:Schur} for the mixed analogue of the Schur-finiteness conjecture.}, denoted by $S(X)$, asserts that the Chow motive $\mathfrak{h}(X)_\bbQ$ is Schur-finite in the sense of Deligne \cite[\S1]{Deligne}.
\end{itemize}
\begin{remark}[Status]
\begin{itemize}
\item[(i)] Thanks to the work of Grothendieck and Kleiman (see \cite{Grothendieck,Kleim1,Kleim}), the conjecture $C^+(X)$ holds when $\mathrm{dim}(X) \leq 2$, and also for abelian varieties. Moreover, this conjecture is stable under products.
\item[(ii)] Thanks to the work of Lieberman \cite{Lieberman}, the conjecture $D(X)$ holds when $\mathrm{dim}(X)\leq 4$, and also for abelian varieties.
\item[(iii)] Thanks to the work Voevodsky \cite{Voevodsky-IMRN} and Voisin \cite{Voisin-nilpotence}, the conjecture $V(X)$ holds when $\mathrm{dim}(X)\leq 2$. Thanks to the work of Kahn-Sebastian \cite{KS}, the conjecture $V(X)$ holds moreover when $X$ is an abelian $3$-fold.
\item[(iv)] Thanks to the work of Kimura \cite{Kimura} and Shermenev \cite{Shermenev}, the conjecture $S(X)$ holds when $\mathrm{dim}(X)\leq 1$, and also for abelian varieties.
\end{itemize}
\end{remark}
Let $k=\bbF_q$ be a finite base field of characteristic $p>0$. Given a smooth proper $k$-scheme $X$ and a prime number $l\neq p$, recall from \cite{Tate-motives, Tate} that the Tate conjecture, denoted by $T^l(X)$, asserts that the cycle class map is surjective:
$$
\cZ^\ast(X)_{\bbQ_l} \too H^{2\ast}_{l\text{-}\text{adic}}(X_{\overline{k}}, \bbQ_l(\ast))^{\mathrm{Gal}(\overline{k}/k)}\,.
$$
\begin{remark}[Status]
Thanks to the work of Tate \cite{Tate}, the conjecture $T^l(X)$ holds when $\mathrm{dim}(X)\leq 1$, and also for abelian varieties. Thanks to the work of several other people (consult Totaro's survey \cite{Totaro}), the conjecture $T^l(X)$ holds moreover when $X$ is a $K3$-surface (and $p\neq 2$).
\end{remark}
\subsection{Noncommutative counterparts}\label{sub:counterparts}
Let $k$ be a base field of characteristic zero. Recall from \S\ref{sec:pure} that periodic cyclic homology descends to the category of noncommutative Chow motives yielding a functor $HP^\pm\colon \NChow(k)_\bbQ \to \mathrm{Vect}_{\bbZ/2}(k)$. Given a smooth proper dg category $\cA$, consider the following $\bbQ$-vector spaces
$$
K_0(\cA)_\bbQ/_{\!\sim ?}:= \Hom_{?}(U(k)_\bbQ, U(\cA)_\bbQ)\,,
$$
where $?$ belongs to $\{\mathrm{nil}, \mathrm{hom}, \mathrm{num}\}$ and $\{\NVoev(k)_\bbQ, \NHom(k)_\bbQ, \NNum(k)_\bbQ\}$, respectively. Under these notations, the important conjectures in \S\ref{sub:recollections2} admit the following noncommutative counterparts:

\vspace{0.1cm}

{\bf Conjecture $C^+_{\mathrm{nc}}(\cA)$:} The even K\"unneth projector $\pi^+_\cA$ of $HP^\pm(\cA)$ is {\em algebraic}, \ie there exists an endomorphism $\underline{\pi}^+_\cA$ of $U(\cA)_\bbQ$ such that $HP^\pm(\underline{\pi}^+_\cA)=\pi^+_\cA$.

\vspace{0.1cm}

{\bf Conjecture $D_{\mathrm{nc}}(\cA)$:} The equality $K_0(\cA)_\bbQ/_{\!\sim \mathrm{hom}}=K_0(\cA)_\bbQ/_{\!\sim \mathrm{num}}$ holds.

\vspace{0.1cm}

{\bf Conjecture $V_{\mathrm{nc}}(\cA)$:} The equality $K_0(\cA)_\bbQ/_{\!\sim \mathrm{nil}}=K_0(\cA)_\bbQ/_{\!\sim \mathrm{num}}$ holds.

\vspace{0.1cm}

{\bf Conjecture $S_{\mathrm{nc}}(\cA)$:} The noncommutative Chow motive $U(\cA)_\bbQ$ is Schur-finite.

\vspace{0.1cm}

Let $k=\bbF_q$ be a finite base field of characteristic $p>0$. Given a smooth proper dg category $\cA$ and a prime number $l\neq p$, consider the following abelian groups 
\begin{eqnarray}\label{eq:abeliangroups}
\Hom\left(\bbZ(l^\infty), \pi_{-1} L_{KU}K(\cA\otimes_{\bbF_q} \bbF_{q^n})\right) && n \geq 1\,,
\end{eqnarray}
where $\bbZ(l^\infty)$ stands for the Pr\"ufer $l$-group and $L_{KU}K(\cA\otimes_k k_n)$ for the Bousfield localization of the algebraic $K$-theory spectrum $K(\cA\otimes_{\bbF_q} \bbF_{q^n})$ with respect to topological complex $K$-theory $KU$. Under these notations, Tate's conjecture admits the following noncommutative counterpart:

\vspace{0.1cm}

{\bf Conjecture $T^l_{\mathrm{nc}}(\cA)$:} The abelian groups \eqref{eq:abeliangroups} are zero.

\vspace{0.1cm}

We now relate the conjectures in \S\ref{sub:recollections2} with their noncommutative counterparts:
\begin{theorem}\label{thm:conjectures}
Given a smooth proper $k$-scheme $X$, we have the equivalences:
\begin{eqnarray}
C^+(X) & \Leftrightarrow & C^+_{\mathrm{nc}}(\perf_\dg(X)) \label{eq:equiv1} \\
D(X) & \Leftrightarrow & D_{\mathrm{nc}}(\perf_\dg(X)) \label{eq:equiv2} \\
V(X) & \Leftrightarrow & V_{\mathrm{nc}}(\perf_\dg(X)) \label{eq:equiv3} \\
S(X) & \Leftrightarrow & S_{\mathrm{nc}}(\perf_\dg(X))\label{eq:equiv4} \\
T^l(X) & \Leftrightarrow & T^l_{\mathrm{nc}}(\perf_\dg(X))\label{eq:equiv5}\,.
\end{eqnarray}
\end{theorem}
Morally speaking, Theorem \ref{thm:conjectures} shows that the important conjectures in \S\ref{sub:recollections2} belong not only to the realm of algebraic geometry but also to the broad noncommutative setting of smooth proper dg categories. Consult \cite[\S5]{book}, and the references therein, for the implications $\Rightarrow$ in \eqref{eq:equiv1}-\eqref{eq:equiv2} and also for the equivalences \eqref{eq:equiv3}-\eqref{eq:equiv4}. The converse implications $\Leftarrow$ in \eqref{eq:equiv1}-\eqref{eq:equiv2} were established in \cite[Thm.~1.1]{note-CD}. Finally, the equivalence \eqref{eq:equiv5} was proved in \cite[Thm.~1.2]{Tate_Tabuada}.
\subsubsection{Homological projective duality}\label{sub:HPD}
For a survey on homological projective duality (=HPD), we invite the reader to consult Kuznetsov's ICM address \cite{ICM-Kuznetsov}. Let $X$ be a smooth projective $k$-scheme equipped with a line bundle $\cL_X(1)$; we write $X \to \bbP(W)$ for the associated morphism where $W:=H^0(X,\cL_X(1))^\ast$. Assume that the category $\perf(X)$ admits a Lefschetz decomposition $\langle \bbA_0, \bbA_1(1), \ldots, \bbA_{i-1}(i-1)\rangle$ with respect to $\cL_X(1)$ in the sense of \cite[Def.~4.1]{KuznetsovHPD}. Following \cite[Def.~6.1]{KuznetsovHPD}, let $Y$ be the HP-dual of $X$, $\cL_Y(1)$ the HP-dual line bundle, and $Y\to \bbP(W^\ast)$ the morphism associated to $\cL_Y(1)$. Given a linear subspace $L \subset W^\ast$, consider the linear sections $X_L:=X\times_{\bbP(W^\ast)} \bbP(L)$ and $Y_L:=Y \times_{\bbP(W)} \bbP(L^\perp)$. The next result, whose proof makes use of Theorem \ref{thm:conjectures}, is obtained by concatenating \cite[\S5.3-5.4]{book} with \cite[Thm.~1.4]{note-CD}\cite[Thm.~1.1]{Schur}\cite[Thm.~1.3]{Tate_Tabuada}:
\begin{theorem}[HPD-invariance\footnote{Consult Theorem \ref{thm:HPD2} for another HPD-invariance type result.}]\label{thm:HPD}
Let $X$ and $Y$ be as above. Assume that $X_L$ and $Y_L$ are smooth, that $\mathrm{dim}(X_L)=\mathrm{dim}(X) -\mathrm{dim}(L)$, that $\mathrm{dim}(Y_L)=\mathrm{dim}(Y)- \mathrm{dim}(L^\perp)$, and that the following conjectures hold
\begin{equation}\label{eq:conjectures}
C^+_{\mathrm{nc}}(\bbA_{0, \dg}) \quad \quad D_{\mathrm{nc}}(\bbA_{0, \dg}) \quad \quad V_{\mathrm{nc}}(\bbA_{0, \dg})\quad \quad S_{\mathrm{nc}}(\bbA_{0, \dg})\quad \quad T^l_{\mathrm{nc}}(\bbA_{0, \dg})\,,
\end{equation}
where $\bbA_{0,\dg}$ stands for the dg enhancement of $\bbA_0$ induced by $\perf_\dg(X)$. Under these assumptions, we have the following equivalences of conjectures:
\begin{eqnarray*}
?(X_L) \Leftrightarrow \,\,?(Y_L) &\text{with}& ?\in \{C^+, D, V, S, T^l\}\,.
\end{eqnarray*}
\end{theorem}
\begin{remark}
The conjectures \eqref{eq:conjectures} hold, for example, whenever the triangulated category $\bbA_0$ admits a full exceptional collection (this is the case in all the examples in the literature). Furthermore, Theorem \ref{thm:HPD} holds more generally when $Y$ (or $X$) is singular. In this case, we need to replace $Y$ by a noncommutative resolution of singularities in the sense of \cite[\S2.4]{ICM-Kuznetsov}.
\end{remark}
Theorem \ref{thm:HPD} shows that the conjectures in \S\ref{sub:recollections2} are invariant under homological projective duality. As a consequence, we obtain the following practical result:
\begin{corollary}\label{cor:HPD}
Let $X_L$ and $Y_L$ be smooth linear sections as in Theorem \ref{thm:HPD}.
\begin{itemize}
\item[(a)] If $\mathrm{dim}(Y_L)\leq 2$, then the conjectures $C^+(X_L)$ and $V(X_L)$ hold.
\item[(b)] If $\mathrm{dim}(Y_L)\leq 4$, then the conjecture $D(X_L)$ holds.
\item[(c)] If $\mathrm{dim}(Y_L)\leq 1$, then the conjectures $S(X_L)$ and $T^l(X_L)$ hold.
\end{itemize}
\end{corollary}
By applying Corollary \ref{cor:HPD} to the Veronese-Clifford duality, to the spinor duality, to the Grassmannian-Pfaffian duality, to the determinantal duality, and to other (incomplete) HP-dualities (see \cite[\S4]{ICM-Kuznetsov}), we obtain a proof of the  conjectures in \S\ref{sub:recollections2} in several new cases; consult \cite{Crelle, note-CD,Schur,Tate_Tabuada} for details. In the particular case of the Veronese-Clifford duality, Corollary \ref{cor:HPD} leads furthermore to an alternative proof of the Tate conjecture for smooth complete intersections of two quadrics (the original (geometric) proof, based on the notion of variety of maximal planes, is due to Reid \cite{Reid}); consult \cite[Thm.~1.7]{Tate_Tabuada} for details.
\subsubsection{Algebraic stacks}\label{sub:stacks}
Theorem \ref{thm:conjectures} allows us to easily extend the important conjectures in \S\ref{sub:recollections2} from smooth proper schemes to smooth proper algebraic stacks $\cX$ by setting $?(\cX):=?_{\mathrm{nc}}(\perf_\dg(\cX))$, where $?\in \{C^+, D, V, S, T^l\}$. The next result, obtained by combining \cite[Thm.~9.2]{Orbifold} with \cite[Thm.~1.9]{Tate_Tabuada}, proves these conjectures in the case of ``low-dimensional'' orbifolds; consult \cite{note-CD} for further examples of algebraic stacks satisfying these conjectures.
\begin{theorem}\label{thm:conj-orbifold} Let $G$ be a finite group, $X$ a smooth projective $k$-scheme equipped with a $G$-action, and $\cX:=[X/G]$ the associated orbifold.
\begin{itemize}
\item[(a)] The conjectures $C^+(\cX)$ and $V(\cX)$ hold when $\mathrm{dim}(X)\leq 2$. The conjecture $C^+(\cX)$ also holds when $G$ acts by group homomorphisms on an abelian variety.
\item[(b)] The conjecture $D(\cX)$ holds when $\mathrm{dim}(X)\leq 4$.
\item[(c)] The conjectures $S(\cX)$ and $T^l(\cX)$ hold when $\mathrm{dim}(X)\leq 1$.
\end{itemize}
\end{theorem}
Roughly speaking, Theorem \ref{thm:conj-orbifold} shows that the above conjectures are ``insensitive'' to the $G$-action. Among other ingredients, its proof makes use~of~Theorem~\ref{thm:orbifold}.
\begin{remark}[Generalizations]
Theorem \ref{thm:conj-orbifold} holds more generally under the assumption that the conjectures in \S\ref{sub:recollections2} are satisfied by the fixed point locus $\{X^\sigma\}_\sigma$, with $\sigma \in \vphi$. For example, the conjecture $T^l(\cX)$ also holds when $X$ is an abelian surface and the group $G=C_2$ acts by the involution $a \mapsto -a$.
\end{remark}
\section{Noncommutative motivic Galois groups}\label{sec:Galois}
Let $F$ be a field of characteristic zero and $\NNum(k)_F$ the abelian category of numerical motives. The next result was proved in \cite[Thm.~6.4]{book} and \cite[Thm.~7.1]{positive}:
\begin{theorem}\label{thm:super}
The category $\NNum(k)_F$ is super-Tannakian in the sense of Deligne \cite{Deligne}. When $F$ is algebraically closed, $\NNum(k)_F$ is neutral super-Tannakian.
\end{theorem}
By combining Theorem \ref{thm:super} with Deligne's super-Tannakian formalism \cite{Deligne}, we obtain an affine super-group $F$-scheme $\mathrm{sGal}(\NNum(k)_F)$ called the {\em noncommutative motivic Galois super-group}. The following result relates this super-group with the (classical) motivic Galois super-group $\mathrm{sGal}(\Num(k)_F)$:
\begin{theorem}[{\cite[Thm.~7.4]{positive}}]\label{thm:Galois}
Assume that $F$ is algebraically closed. Then, there exists a faithfully flat morphism of affine super-group $F$-schemes
$$ \mathrm{sGal}(\NNum(k)_F) \twoheadrightarrow \mathrm{Ker} (\mathrm{sGal}(\Num(k)_F) \stackrel{t^\ast}{\twoheadrightarrow} \bbG_m)\,,$$
where $\bbG_m$ stands for the multiplicative (super-)group scheme and $t$ for the inclusion of the category of Tate motives into numerical motives.
\end{theorem}
Theorem \ref{thm:Galois} was envisioned by Kontsevich; see his seminal talk \cite{IAS}. Intuitively speaking, it shows that the ``$\otimes$-symmetries'' of the commutative world which can be lifted to the noncommutative world are precisely those which become trivial when restricted to Tate motives. Theorem \ref{thm:Galois} also holds when $F$ is {\em not} algebraically closed. However, in this case the super-group schemes are only defined over a (very big) commutative $F$-algebra.
\begin{remark}[Simplification]
The analogue of Theorem \ref{thm:Galois}, with $k$ of characteristic zero, was proved in \cite[Thm.~6.7(ii)]{book}. However, therein we assumed the noncommutative counterparts of the standard conjectures of type $C^+$ and $D$ and moreover used Deligne-Milne's theory of Tate-triples. In contrast, Theorem \ref{thm:Galois} is unconditional and its proof avoids the use of Tate-triples; consult \cite[\S7]{positive} for details.
\end{remark}
\subsection*{Base-change}
Recall from \cite[\S6]{book} the definition of the (conditional\footnote{We assume the noncommutative counterparts of the standard conjectures of type $C^+$ and $D$.}) noncommutative motivic Galois group $\mathrm{Gal}(\NNum^\dagger(k)_F)$.
\begin{theorem}[{\cite[Thm.~2.2]{rigidity}}]\label{thm:change}
Given a primary field extension $l/k$, the induced base-change functor $-\otimes_k l \colon \NNum^\dagger(k)_F \to \NNum^\dagger(l)_F$ gives rise to a faithfully flat morphism of affine group $F$-schemes $\mathrm{Gal}(\NNum^\dagger(l)_F) \to \mathrm{Gal}(\NNum^\dagger(k)_F)$.
\end{theorem}
Roughly speaking, Theorem \ref{thm:Galois} shows that every ``$\otimes$-symmetry'' of the category of noncommutative numerical $k$-linear motives can be extended to a ``$\otimes$-symmetry'' of the category of noncommutative $l$-linear motives. Among other ingredients, its proof makes use of Theorem \ref{thm:rigidity}. In the particular case of an extension of algebraically closed fields $l/k$, the commutative counterpart of Theorem \ref{thm:change} was established by Deligne-Milne in \cite[Prop.~6.22(b)]{DM}.
\section{Localizing invariants}\label{sec:localizing}
Recall from \cite[\S8.1]{book} the notion of a {\em short exact sequence} of dg categories in the sense of Drinfeld/Keller. In \S\ref{sub:ses} we describe a key structural property of these short exact sequences and explain its implications to secondary $K$-theory.

Recall from \cite[\S8.5.1]{book} the construction of the {\em universal localizing $\bbA^1$-homotopy invariant} of dg categories $\mathrm{U}\colon \dgcat(k) \to \Nmot(k)$; in {\em loc. cit.} we used the explicit notation $\mathrm{U}^{\bbA^1}_\loc\colon \dgcat(k) \to \Nmot_\loc^{\bbA^1}(k)$. In \cite[\S8.5.3]{book} we described the behavior of $\mathrm{U}$ with respect to dg orbit categories and dg cluster categories. In \S\ref{sub:Gysin}-\ref{sub:NCproj} we describe the behavior of $\mathrm{U}$ with respect to open/closed scheme decompositions, corner skew Laurent polynomial algebras, and noncommutative projective schemes. As explained in \cite[Thm.~8.25]{book}, all the results in \S\ref{sub:Gysin}-\ref{sub:NCproj} are {\em motivic} in the sense that they hold similarly for every localizing $\bbA^1$-homotopy invariant such as mod$\text{-}n$ algebraic $K$-theory (when $1/n \in k$), homotopy $K$-theory, \'etale $K$-theory, periodic cyclic homology\footnote{Periodic cyclic homology is not a localizing $\bbA^1$-homotopy invariant in the sense of \cite[\S8.5]{book} because it does not preserves filtered (homotopy) colimits. Nevertheless, all the results of \S\ref{sub:Gysin}-\ref{sub:NCproj} hold similarly for periodic cyclic homology.} (when $\mathrm{char}(k)=0$),~etc. The results of \S\ref{sub:NCproj} do {\em not} require $\bbA^1$-homotopy invariance and so they hold for every localizing invariant; see~Remark~\ref{rk:localizing}.

\begin{notation} Given a $k$-scheme $X$ (or more generally an algebraic stack $\cX$), we will write $\mathrm{U}(X)$ instead of $\mathrm{U}(\perf_\dg(X))$.
\end{notation}
\subsection{Short exact sequences}\label{sub:ses}
Recall from \cite[\S8.4]{book} the notion of a {\em split short exact sequence} of dg categories $0 \to \cA\to \cB \to \cC \to 0$. Up to Morita equivalence, this data is equivalent to inclusions of dg categories $\cA , \cC \subseteq \cB$ yielding a semi-orthogonal decomposition of triangulated categories $\dgHo(\cB)=\langle \dgHo(\cA), \dgHo(\cC)\rangle$ in the sense of Bondal-Orlov \cite{BO}; by definition, the category $\dgHo(\cA)$ has the same objects as $\cA$ and morphisms $\dgHo(\cA)(x,y):=H^0(\cA(x,y))$.
\begin{theorem}[{\cite[Thm.~4.4]{Secondary}}]\label{thm:split}
Let $0 \to \cA \to \cB \to \cC \to 0$ be a short exact sequence of dg categories in the sense of Drinfeld/Keller. If $\cA$ is smooth and proper and $\cB$ is proper, then the short exact sequence is split.
\end{theorem}
Morally speaking, Theorem \ref{thm:split} shows that the smooth proper dg categories behave as ``injective'' objects. In the setting of triangulated categories, this conceptual idea goes back to the pioneering work of Bondal-Kapranov \cite{BK}.
\subsubsection{Secondary $K$-theory}\label{sub:secondary}
Two decades ago, Bondal-Larsen-Lunts introduced in \cite{BLL} the {\em Grothendieck ring of smooth proper dg categories $\cP\cT(k)$}. This ring is defined by generators and relations. The generators are the Morita equivalence classes of smooth proper dg categories\footnote{Bondal-Larsen-Lunts worked originally with (pretriangulated) dg categories. In this generality, the classical Eilenberg's swindle argument implies that the Grothendieck ring is trivial.} and the relations $[\cB]=[\cA]+[\cC]$ arise from semi-orthogonal decompositions $\dgHo(\cB)=\langle \dgHo(\cA), \dgHo(\cC) \rangle$. The multiplication law is induced by the tensor product of dg categories. One decade ago, To\"en introduced in \cite{Toen} a ``categorified'' version of the Grothendieck ring named {\em secondary Grothendieck ring $K_0^{(2)}(k)$}. By definition, $K_0^{(2)}(k)$ is the~quotient of the free abelian group on the Morita equivalence classes of smooth proper dg categories by the relations $[\cB]=[\cA]+[\cC]$ arising from short exact sequences $0 \to \cA \to \cB \to \cC\to 0$. The multiplication law is also induced by the tensor product of dg categories.

Theorem \ref{thm:split} directly leads to the following result:
\begin{corollary}\label{cor:secondary}
The rings $\cP\cT(k)$ and $K_0^{(2)}(k)$ are isomorphic. 
\end{corollary}
Morally speaking, Corollary \ref{cor:secondary} shows that the secondary Grothendieck ring is {\em not} a new mathematical notion. 

By construction, the universal additive invariant $U$ (see \S\ref{sec:additive}) sends semi-orthogonal decompositions to direct sums. Therefore, it gives rise to a ring homomorphism $\cP\cT(k)\to K_0(\NChow(k))$. Making use of Corollary \ref{cor:new}, we then obtain the result:
\begin{corollary}\label{cor:secondary1}
The following map is injective:
\begin{eqnarray}\label{eq:map}
\mathrm{Br}(k) \too \cP\cT(k)\simeq K_0^{(2)}(k) && [A] \mapsto [A]\,.
\end{eqnarray}
\end{corollary}
The map \eqref{eq:map} may be understood as the ``categorification'' of the canonical map from the Picard group $\mathrm{Pic}(k)$ to the Grothendieck ring $K_0(k)$. In contrast with $\mathrm{Pic}(k) \to K_0(k)$, the map \eqref{eq:map} does not seems to admit a ``determinant'' map in the converse direction. Nevertheless, Corollary \ref{cor:secondary1} shows that this~map~is~injective.
\begin{remark}[Generalizations]
Similarly to Remark \ref{rk:generalization1}, Corollaries \ref{cor:secondary}-\ref{cor:secondary1} hold more generally with $k$ replaced by a base $k$-scheme $X$.
\end{remark}
\subsection{Gysin triangle}\label{sub:Gysin}
Let $X$ be a smooth $k$-scheme, $i\colon Z \hookrightarrow X$ a smooth closed subscheme, and $j\colon V \hookrightarrow X$ the open complement of $Z$.
\begin{theorem}[{\cite[Thm.~1.9]{Gysin}}]\label{thm:Gysin}
We have an induced distinguished ``Gysin'' triangle
\begin{equation}\label{eq:Gysin-mot1}
\mathrm{U}(Z) \stackrel{\mathrm{U}(i_\ast)}{\too} \mathrm{U}(X) \stackrel{\mathrm{U}(j^\ast)}{\too} \mathrm{U}(V) \stackrel{\partial}{\too} \mathrm{U}(Z)[1]\,,
\end{equation}
where $i_\ast$, resp. $j^\ast$, stands for the push-forward, resp. pull-back, dg functor.
\end{theorem}
\begin{remark}[Generalizations]
As explained in \cite[\S7]{Gysin}, Theorem \ref{thm:Gysin} holds not only for smooth schemes but also for smooth algebraic spaces in the sense of Artin.
\end{remark}
Roughly speaking, Theorem \ref{thm:Gysin} shows that the difference between the localizing $\bbA^1$-homotopy invariants of $X$ and of $V$ is completely determined by the closed subscheme $Z$. Consult Remark \ref{rk:motivic1}, resp. \ref{rk:motivic2}, for the relation between \eqref{eq:Gysin-mot1} and the motivic Gysin triangles constructed by Morel-Voevodsky, resp. Voevodsky.
\subsubsection{Quillen's localization theorem}
Homotopy $K$-theory is a localizing $\bbA^1$-homotopy invariant which agrees with Quillen's algebraic $K$-theory when restricted to smooth $k$-schemes. Therefore, Theorem \ref{thm:Gysin} leads to the $K$-theoretical localization theorem
\begin{equation}\label{eq:localization}
K(Z) \stackrel{K(i_\ast)}\too K(X) \stackrel{K(j^\ast)}{\too} K(V) \stackrel{\partial}{\too} K(Z)[1]
\end{equation}
originally established by Quillen in \cite[Chapter 7 \S3]{Quillen}. Among other ingredients, Quillen's proof makes essential use of 
d\'evissage. The proof of Theorem \ref{thm:Gysin}, and hence of \eqref{eq:localization}, is quite different and avoids the use of d\'evissage.
\subsubsection{Six-term exact sequence in de Rham cohomology}\label{sub:deRham}
Periodic cyclic homology is a localizing $\bbA^1$-homotopy invariant (when $\mathrm{char}(k)=0$). Moreover, thanks to the Hochschild-Kostant-Rosenberg theorem, we have an isomorphism of $\bbZ/2$-graded $k$-vector spaces $HP^\pm(X) \simeq (\bigoplus_{n\,\mathrm{even}} H^n_{dR}(X), \bigoplus_{n\,\mathrm{odd}}H^n_{dR}(X))$, where $H^\ast_{dR}$ stands for de Rham cohomology. Furthermore, the maps $i$ and $j$ give rise to homomorphisms $H^n_{dR}(i_\ast)\colon H^n_{dR}(Z) \to H^{n+2c}_{dR}(X) $ and $H^n_{dR}(j^\ast)\colon H^n_{dR}(X) \to H^n_{dR}(V)$, where $c$ stands for the codimension of $i$. Therefore, Theorem \ref{thm:Gysin} leads to the following six-term exact sequence in de Rham cohomology:
$$
\xymatrix{
\bigoplus_{n\,\mathrm{even}}H^n_{dR}(Z) \ar[rr]^-{\bigoplus_nH^n_{dR}(i_\ast)} && \bigoplus_{n\,\mathrm{even}}H^n_{dR}(X) \ar[rr]^-{\bigoplus_nH^n_{dR}(j^\ast)} && \bigoplus_{n\,\,\mathrm{even}}H^n_{dR}(V) \ar[d]^-\partial \\
\bigoplus_{n\,\mathrm{odd}}H^n_{dR}(V) \ar[u]^-\partial && \bigoplus_{n\,\mathrm{odd}}H^n_{dR}(X) \ar[ll]^-{\bigoplus_n H^n_{dR}(j^\ast)} && \bigoplus_{n\,\,\mathrm{odd}}H^n_{dR}(Z) \ar[ll]^-{\bigoplus_nH^n_{dR}(i_\ast)} \,.  
}
$$
This exact sequence is the ``$2$-periodization'' of the
Gysin long exact sequence on de
Rham~cohomology~originally~constructed~by~Hartshorne~in~\cite[Chapter
II~\S3]{Hartshorne}.  
\subsubsection{Reduction to projective schemes}
As a byproduct of Theorem \ref{thm:Gysin}, the study of the localizing $\bbA^1$-homotopy invariants of smooth $k$-schemes can be reduced to the study of the localizing $\bbA^1$-homotopy invariants of smooth {\em projective} $k$-schemes:
\begin{theorem}[{\cite[Thm.~2.1]{Gysin}}]\label{thm:reduction}
Let $X$ a smooth $k$-scheme.
\begin{itemize}
\item[(i)] If $\mathrm{char}(k)=0$, then $\mathrm{U}(X)$ belongs to the smallest triangulated subcategory of $\Nmot(k)$ containing the objects $\mathrm{U}(Y)$, with $Y$ a smooth projective $k$-scheme.
\item[(ii)] If $k$ is a perfect field of characteristic $p>0$, then $\mathrm{U}(X)_{1/p}$ belongs to the smallest thick triangulated subcategory of $\Nmot(k)_{1/p}$ containing the objects $\mathrm{U}(Y)_{1/p}$, with $Y$ a smooth projective $k$-scheme.
\end{itemize}
\end{theorem}
Among other ingredients, the proof of item (i), resp. item (ii), of Theorem \ref{thm:reduction} makes use of resolution of singularities, resp. of Gabber's refined version of de Jong's theory of alterations; consult \cite[\S8]{Gysin} for details.
\begin{remark}[Dualizable objects]\label{rk:dualizable}
Given a smooth projective $k$-scheme $Y$, the associated dg category $\perf_\dg(Y)$ is smooth and proper; see \cite[Example~1.42]{book}. Therefore, since the universal localizing $\bbA^1$-homotopy invariant $\mathrm{U}$ is symmetric monoidal, it follows from \cite[Thm.~1.43]{book} that $\mathrm{U}(Y)$ is a dualizable object of the symmetric monoidal category $\Nmot(k)$. Given a smooth $k$-scheme $X$, the associated dg category $\perf_\dg(X)$ is smooth but {\em not} necessarily proper. Nevertheless, Theorem \ref{thm:reduction} implies that $\mathrm{U}(X)$, resp. $\mathrm{U}(X)_{1/n}$, is still a dualizable object of the symmetric monoidal category $\Nmot(k)$, resp. $\Nmot(k)_{1/p}$.
\end{remark}
\subsection{Corner skew Laurent polynomial algebras}\label{sub:corner}
Let $A$ be a unital $k$-algebra, $e$ an idempotent of $A$, and $\phi\colon A\stackrel{\sim}{\to} eAe$ a ``corner'' isomorphism. The associated {\em corner skew Laurent polynomial algebra $A[t_+,t_-;\phi]$} is defined as follows: the elements are formal expressions $t^m_- a_{-m} + \cdots + t_- a_{-1} + a_0 + a_1 t_+ \cdots + a_n t_+^n$
with $a_{-i} \in \phi^i(1)A$ and $a_i \in A\phi^i(1)$ for every $i \geq 0$; the addition is defined componentwise; the multiplication is determined by the distributive law and by the relations $t_-t_+=1$, $t_+t_-=e$, $at_-=t_-\phi(a)$ for every $a \in A$, and $t_+a = \phi(a)t_+$ for every $a \in A$. Note that $A[t_+,t_-;\phi]$ admits a canonical $\bbZ$-grading with $\mathrm{deg}(t_\pm)=\pm 1$. As proved in \cite[Lem.~2.4]{Fractional}, the corner skew Laurent polynomial algebras can be characterized as those $\bbZ$-graded algebras $C=\bigoplus_{n \in \bbZ}C_n$ containing elements $t_+\in C_1$ and $t_-\in C_{-1}$ such that $t_-t_+=1$. Concretely, we have $C=A[t_+,t_-;\phi]$ with $A:=C_0$, $e:=t_+t_-$, and $\phi\colon C_0 \to t_+ t_- C_0 t_+ t_-$ given by $c_0\mapsto t_+ c_0 t_-$.
\begin{example}[Skew Laurent polynomial algebras]
When $e=1$, $A[t_+,t_-;\phi]$ reduces to the classical skew Laurent polynomial algebra $A \rtimes_\phi \bbZ$. In the particular case where $\phi$ is the identity, $A\rtimes_\phi \bbZ$ reduces furthermore to $A[t,t^{-1}]$.
\end{example}
\begin{example}[Leavitt algebras]
Following \cite{Leavitt}, the {\em Leavitt algebra $L_n$, $n\geq 0$,} is the $k$-algebra generated by elements $x_0, \ldots, x_n, y_0, \ldots, y_n$ subject to the relations $y_i x_j = \delta_{ij}$ and $\sum^n_{i=0} x_i y_i =1$. Note that the canonical $\bbZ$-grading, with $\mathrm{deg}(x_i)=1$ and $\mathrm{deg}(y_i)=-1$, makes $L_n$ into a corner skew Laurent polynomial algebra. Note also that  $L_0\simeq k[t,t^{-1}]$. In the remaining cases $n\geq 1$, $L_n$ is the universal example of a $k$-algebra of {\em module type $(1,n+1)$}, \ie $L_n\simeq L_n^{\oplus (n+1)}$ as right $L_n$-modules.
\end{example}
\begin{example}[Leavitt path algebras]\label{ex:Leavitt-path}
Let $Q=(Q_0,Q_1, s,r)$ be a finite quiver with no sources; $Q_0$ and $Q_1$ stand for the sets of vertices and arrows, respectively, and $s$ and $r$ for the source and target maps, respectively. Consider the double quiver $\overline{Q}=(Q_0,Q_1\cup Q_1^\ast, s,r)$ obtained from $Q$ by adding an arrow $\alpha^\ast$, in the converse direction, for each arrow $\alpha \in Q_1$. The {\em Leavitt path algebra $L_Q$ of $Q$} is the quotient of the quiver algebra $k\overline{Q}$ (which is generated by elements $\alpha \in Q_1\cup Q_1^\ast$ and $e_i$ with $i \in Q_0$) by the Cuntz-Krieger's relations: $\alpha^\ast \beta = \delta_{\alpha \beta} e_{r(\alpha)}$ for every $\alpha, \beta \in Q_1$ and $\sum_{\{\alpha \in Q_1|s(\alpha)=i\}} \alpha \alpha^\ast =e_i$ for every non-sink $i \in Q_0$. Note that $L_Q$ admits a canonical $\bbZ$-grading with $\mathrm{deg}(\alpha)=1$ and $\mathrm{deg}(\alpha^\ast)=-1$. For every vertex $i \in Q_0$ choose an arrow $\alpha_i$ such that $r(\alpha_i)=i$ and consider the associated elements $t_+:=\sum_{i \in Q_0}\alpha_i$ and $t_-:=t_+^\ast$. Since $\mathrm{deg}(t_\pm)=\pm1$ and $t_-t_+=1$, $L_Q$ is an example of a corner skew Laurent polynomial algebra. In the particular case where $Q$ is the quiver with one vertex and $n+1$ arrows, $L_Q$ reduces to the Leavitt algebra $L_n$. 
\end{example}
\begin{theorem}[{\cite[Thm.~3.1]{Corner}}]\label{thm:corner}
We have an induced distinguished triangle
\begin{equation*}
\mathrm{U}(A) \stackrel{\id - \mathrm{U}({}_\phi A)}{\too} \mathrm{U}(A) \too \mathrm{U}(A[t_+, t_-; \phi]) \stackrel{\partial}{\too} \mathrm{U}(A)[1]\,,
\end{equation*}
where ${}_\phi A$ stands for the $A\text{-}A$-bimodule associated to $\phi$.
\end{theorem}
Roughly speaking, Theorem \ref{thm:corner} shows that $A[t_+, t_-; \phi]$ may be considered as a model for the orbits of the $\bbN$-action on $\mathrm{U}(A)$ induced by the endomorphism $\mathrm{U}({}_\phi A)$.
\subsubsection{Leavitt path algebras}
Let $Q=(Q_0,Q_1,s,r)$ be a quiver as in Example \ref{ex:Leavitt-path}, with $v$ vertices and $v'$ sinks. Assume that the set $Q_0$ is ordered with the first $v'$ elements corresponding to the sinks. Let $I'_Q$ be the incidence matrix of $Q$, $I_Q$ the matrix obtained from $I'_Q$ by removing the first $v'$ rows (which are zero), and $I_Q^t$ the transpose of $I_Q$. Under these notations, Theorem \ref{thm:corner} (concerning the Leavitt path algebra $L_Q$) admits the following refinement:
\begin{theorem}[{\cite[Thm.~3.7]{Corner}}]\label{thm:corner2}
We have an induced distinguished triangle:
\begin{equation}\label{eq:triangle2}
\mathrm{U}(k)^{(v-v') \oplus} \stackrel{\binom{0}{\id} - I^t_Q}{\too} \mathrm{U}(k)^{v \oplus} \too \mathrm{U}(L_Q) \stackrel{\partial}{\too} \mathrm{U}(k)[1]^{(v-v') \oplus}\,.
\end{equation}
\end{theorem}
Roughly speaking, Theorem \ref{thm:corner2} shows that all the information about the localizing $\bbA^1$-homotopy invariants of Leavitt path algebras $L_Q$ is encoded in the incidence matrix of the quiver $Q$. As an application, Theorem \ref{thm:corner2} directly leads to the following explicit model for the mod$\text{-}n$ Moore construction\footnote{Explicit models for the suspension construction, namely the Waldhausen's $S_\bullet$-construction and the Calkin algebra, are described in \cite[\S8.3.2 and \S8.4.4]{book}.}:
\begin{example}[mod$\text{-}n$ Moore construction]
Let $Q$ be the quiver with one vertex and $n+1$ arrows. In this particular case, \eqref{eq:triangle2} reduces to the distinguished triangle
\begin{equation*}
\mathrm{U}(k) \stackrel{n\cdot \id}{\too} \mathrm{U}(k) \too \mathrm{U}(L_n) \stackrel{\partial}{\too} \mathrm{U}(k)[1]\,.
\end{equation*}
This shows that the Leavitt algebra $L_n, n\geq 2$, is a model for the mod$\text{-}n$ Moore object of $\mathrm{U}(k)$. Therefore, since the universal localizing $\bbA^1$-homotopy invariant $\mathrm{U}$ is symmetric monoidal, given a small dg category $\cA$, we conclude that the tensor product $\cA\otimes L_n$ is a model for the mod$\text{-}n$ Moore object of $\mathrm{U}(\cA)$. 
\end{example}
\subsection{Noncommutative projective schemes}\label{sub:NCproj}
Let $A=\bigoplus_{n\geq 0} A_n$ be a $\bbN$-graded Noetherian $k$-algebra. In what follows, we assume that $A$ is {\em connected}, \ie $A_0=k$, and {\em locally finite-dimensional}, \ie $\mathrm{dim}_k(A_n)<\infty$ for every $n$. Following Manin \cite{Manin}, Gabriel \cite{Gabriel}, Artin-Zhang \cite{Artin-Zhang}, and others, the {\em noncommutative projective scheme $\mathrm{qgr}(A)$} associated to $A$ is defined as the quotient abelian category $\mathrm{gr}(A)/\mathrm{tors}(A)$, where $\mathrm{gr}(A)$ stands for the abelian category of finitely generated $\bbZ$-graded (right) $A$-modules and $\mathrm{tors}(A)$ for the Serre subcategory of torsion $A$-modules. This definition was motivated by Serre's celebrated result \cite[Prop.~7.8]{Serre}, which asserts that in the particular case where $A$ is commutative and generated by elements of degree $1$ the quotient category $\mathrm{qgr}(A)$ is equivalent to the abelian category of coherent $\cO_{\mathrm{Proj}(A)}$-modules $\mathrm{coh}(\mathrm{Proj}(A))$. For example, when $A$ is the polynomial $k$-algebra $k[x_1, \ldots, x_d]$, with $\mathrm{deg}(x_i)=1$, we have the equivalence of categories $\mathrm{qgr}(k[x_1, \ldots, x_d])\simeq \mathrm{coh}(\bbP^{d-1})$. For a survey on noncommutative projective geometry, we invite the reader to consult Stafford's ICM address \cite{Stafford}.

Assume that $A$ is Koszul and has finite global dimension $d$. Under these assumptions, the Hilbert series $h_A(t):=\sum_{n\geq 0} \mathrm{dim}_k(A_n)t^n \in \bbZ[\![t]\!]$ is invertible and its inverse $h_A(t)^{-1}$ is a polynomial $1-\beta_1t + \beta_2t^2 - \cdots + (-1)^d\beta_d t^d$ of degree $d$, with $\beta_i$ the dimension of the $k$-vector space $\mathrm{Tor}^A_i(k,k)$ (or $\mathrm{Ext}_A^i(k,k)$).
\begin{example}[Quantum polynomial algebras]
Choose constant elements $q_{ij} \in k^\times$ with $1 \leq i < j \leq d$. Following Manin~\cite[\S1]{Manin-Fourier}, the $\bbN$-graded Noetherian $k$-algebra
$$ A:=k\langle x_1, \ldots, x_d\rangle/\langle x_j x_i - q_{ij} x_i x_j\,|\, 1 \leq i < j \leq d\rangle\,,$$
with $\mathrm{deg}(x_i)=1$, is called the {\em quantum polynomial algebra} associated to $q_{ij}$. This algebra is Koszul,  has global dimension $d$, and $h_A(t)^{-1}=(1-t)^d$.
\end{example}
\begin{example}[Quantum matrix algebras]
Choose a constant element $q \in k^\times$. Following Manin \cite[\S1]{Manin-Fourier}, the $\bbN$-graded Noetherian $k$-algebra $A$ defined as the quotient of $k\langle x_1, x_2, x_3, x_4\rangle$ by the following relations
\begin{eqnarray*}
x_1x_2 = q x_2 x_1 & x_1x_3=qx_3 x_1 &  x_1 x_4 - x_4 x_1 = (q-q^{-1}) x_2 x_3 \\
x_2 x_3 = x_3 x_2 &  x_2 x_4=q x_4 x_2 & x_1 x_4 = q x_4 x_3\,,
\end{eqnarray*}
with $\mathrm{deg}(x_i)=1$, is called the {\em quantum matrix algebra} associated to $q$. This algebra is Koszul, has global dimension $4$, and $h_A(t)^{-1}=(1-t)^4$.
\end{example}
\begin{example}[Sklyanin algebras]
Let $C$ be a smooth elliptic $k$-curve, $\sigma$ an automorphism of $C$ given by translation under the group law, and $\cL$ a line bundle on $C$ of degree $d\geq 3$. We write $\Gamma_\sigma \subset C\times C$ for the graph of $\sigma$ and $W$ for the $d$-dimensional $k$-vector space $H^0(C,\cL)$. Following Feigin-Odesskii \cite{Feigin} and Tate-Van den Bergh \cite[\S1]{TateVdb}, the $\bbN$-graded Noetherian $k$-algebra $A:=T(W)/R$, where
$$R:=H^0(C\times C, (\cL\boxtimes \cL)(-\Gamma_\sigma))\subset H^0(C\times C, \cL \boxtimes \cL)=W\otimes W\,,$$
is called the {\em Sklyanin algebra} associated to the triple $(C,\sigma, \cL)$. This algebra is Koszul, has global dimension $d$, and $h_A(t)^{-1}=(1-t)^d$.
\end{example}
\begin{example}[Homogenized enveloping algebras]
Let $\mathfrak{g}$ be a finite-dimensional Lie algebra. Following Smith \cite[\S12]{Smith}, the $\bbN$-graded Noetherian $k$-algebra
$$ A:= T(\mathfrak{g} \oplus kz)/ \langle\{z\otimes x - x\otimes z \,|\, x \in \mathfrak{g}\}\cup \{x\otimes y - y\otimes x - [x,y]\otimes z \,|\, x, y \in \mathfrak{g}\}\rangle\,,$$
is called the {\em homogenized enveloping algebra} of $\mathfrak{g}$. This algebra is Koszul, has global dimension $d:=\mathrm{dim}(\mathfrak{g})+1$, and $h_A(t)^{-1}=(1-t)^{d}$.
\end{example}
Given a $\bbN$-graded $k$-algebra $A$ as above, let us write $\cD^b_{\mathrm{dg}}(\mathrm{qgr}(A))$ for the canonical dg enhancement of the bounded derived category of $\mathrm{qgr}(A)$. This dg category is, in general, {\em not} proper; see \cite[\S1]{NCProj}. The following result is contained \cite[Thm.~1.2]{NCProj}:
\begin{theorem}\label{thm:NCProj}
We have an induced distinguished triangle
\begin{equation}\label{eq:triangle-NCProj}
\bigoplus^{+\infty}_{-\infty} \mathrm{U}(k) \stackrel{\mathrm{M}}{\too}  \bigoplus^{+\infty}_{-\infty} \mathrm{U}(k) \too \mathrm{U}(\cD^b_{\mathrm{dg}}(\mathrm{qgr}(A))) \stackrel{\partial}{\too} \bigoplus^{+\infty}_{-\infty} \mathrm{U}(k)[1]\,,
\end{equation}
where $\mathrm{M}$ stands for the (infinite) matrix $\mathrm{M}_{ij}:=(-1)^j (-1)^{(i-j)} \beta_{i-j}$. Moreover, when $\beta_d=1$, the triangle \eqref{eq:triangle-NCProj} induces an isomorphism $\mathrm{U}(\cD^b_{\mathrm{dg}}(\mathrm{qgr}(A)))\simeq \mathrm{U}(k)^{\oplus d}$.
\end{theorem}
As proved in \cite[Cor.~0.2]{Zhang}, we have $h_A(t)^{-1}=(1-t)^3$ whenever $d=3$.
\begin{remark}[Localizing invariants]\label{rk:localizing}
The proof of Theorem \ref{thm:NCProj} does {\em not} makes use of $\bbA^1$-homotopy invariance. Consequently, as explained in \cite[Thm.~8.5]{book}, Theorem \ref{thm:NCProj} holds similarly for every localizing invariant in the sense of \cite[Def.~8.3]{book}. Examples of localizing invariants which are {\em not} $\bbA^1$-homotopy invariant include nonconnective algebraic $K$-theory, Hochschild homology, cyclic homology, negative cyclic homology, periodic cyclic homology (when $\mathrm{char}(k)=p>0$), topological Hochschild homology, topological cyclic homology, topological periodic cyclic homology, etc.
\end{remark}
Roughly speaking, Theorem \ref{thm:NCProj} (and Remark \ref{rk:localizing}) shows that the localizing invariants of a noncommutative projective scheme $\mathrm{qgr}(A)$ are completely determined by the Hilbert series $h_A(t)$.
\section{Noncommutative mixed motives}\label{sec:NCmixed}
In this section we assume that the base field $k$ is perfect. Kontsevich introduced in \cite{Miami,finMot,IAS} a certain rigid symmetric monoidal triangulated category of noncommutative mixed motives $\Nmix(k)$. As explained in \cite[\S9.1.1]{book}, this category can be (conceptually) described as the smallest thick triangulated subcategory of $\Nmot(k)$ (see \S\ref{sec:localizing}) containing the objects $\mathrm{U}(\cA)$, with $\cA$ smooth and proper. 

In \S\ref{sub:Picard} we compute the Picard group of the thick triangulated subcategory of $\mathrm{NMix}(k)$ generated by the noncommutative mixed motives of central simple $k$-algebras. Subsections \S\ref{sub:MV}-\ref{sub:Levine} are devoted to the precise relation between the category $\Nmix(k)$ and Morel-Voevodsky's stable $\bbA^1$-homotopy category, Voevodsky's triangulated category of geometric mixed motives, and Levine's triangulated category of mixed motives, respectively. In \S\ref{sec:Schur} we address the Schur-finiteness conjecture in the case of quadric fibrations. Finally, subsection \S\ref{sub:rigidity2} is devoted to the rigidity property of the category of mod-$n$ noncommutative mixed motives.
\subsection{Picard group}\label{sub:Picard}
The computation of the Picard group of the category of noncommutative mixed motives is a major challenge which seems completely out of reach at the present time. However, this major challenge can be met if we restrict ourselves to central simple $k$-algebras. Let $\Nmix_{\mathrm{csa}}(k)$ be the thick triangulated subcategory of $\Nmix(k)$ generated by the noncommutative mixed motives $\mathrm{U}(A)$ of central simple $k$-algebras $A$. Similarly to \S\ref{sub:Brauer}, the equivalence $[A]=[B]\Leftrightarrow \mathrm{U}(A)\simeq \mathrm{U}(B)$ holds for any two central simple $k$-algebras $A$ and $B$. Moreover, following \cite[Thm.~8.28]{book}, we have non-trivial  Ext-groups:
\begin{equation}\label{eq:equiv-last1}
\Hom_{\Nmix(k)}(\mathrm{U}(A), \mathrm{U}(B)[-n]) \simeq K_n(A^\op \otimes B) \quad \quad n\in \bbZ\,.
\end{equation}
This shows that $\Nmix_{\mathrm{csa}}(k)$ contains information not only about the Brauer group $\mathrm{Br}(k)$ but also about all the higher algebraic $K$-theory of central simple $k$-algebras.
\begin{theorem}[{\cite[Thm.~2.22]{Picard}}]\label{thm:Picard}
We have the following isomorphism:
\begin{eqnarray*}
\mathrm{Br}(k) \times \bbZ \stackrel{\sim}{\too} \mathrm{Pic}(\Nmix_{\mathrm{csa}}(k)) && ([A],n) \mapsto \mathrm{U}(A)[n]\,.
\end{eqnarray*}
\end{theorem}
Theorem \ref{thm:Picard} shows that, although $\Nmix_{\mathrm{csa}}(k)$ contains information about all the higher algebraic $K$-theory of central simple $k$-algebras, none of the noncommutative mixed motives which are built using the non-trivial Ext-groups \eqref{eq:equiv-last1} is $\otimes$-invertible.
\subsection{Morel-Voevodsky's motivic category}\label{sub:MV}
Morel-Voevodsky introduced in \cite{MV,Voevodsky-ICM} the stable $\bbA^1$-homotopy category of $(\bbP^1,\infty)$-spectra $\mathrm{SH}(k)$. By construction, we have a symmetric monoidal functor $\Sigma^\infty(-_+)\colon \mathrm{Sm}(k) \to \mathrm{SH}(k)$ defined on smooth $k$-schemes. Let $\mathrm{KGL} \in \mathrm{SH}(k)$ be the ring $(\bbP^1,\infty)$-spectrum representing homotopy $K$-theory and $\mathrm{Mod}(\mathrm{KGL})$ the homotopy category of $\mathrm{KGL}$-modules.
\begin{theorem}\label{thm:bridge1}
\begin{itemize}
\item[(i)] If $\mathrm{char}(k)=0$, then there exists a fully-faithful, symmetric monoidal, triangulated functor $\Psi$ making the following diagram commute
\begin{equation}\label{eq:diagram-1}
\xymatrix{
\mathrm{Sm}(k) \ar[d]_-{\Sigma^\infty(-_+)} \ar[rrr]^-{X \mapsto \perf_\dg(X)} \ar[drr] &&& \dgcat(k) \ar[d]^-{\mathrm{U}} \\
\mathrm{SH}(k) \ar[d]_-{-\wedge \mathrm{KGL}} && \Nmix(k) \ar[d]_-{(-)^\vee} \ar[r] & \Nmot(k) \ar[d]^-{\uHom(-,\mathrm{U}(k))} \\
\mathrm{Mod}(\mathrm{KGL}) \ar[rr]_-\Psi && \Nmix(k)^\oplus \ar[r] & \Nmot(k) \,,
}
\end{equation}
where $\uHom(-,-)$ stands for the internal-Hom of the closed symmetric monoidal category $\Nmot(k)$, $(-)^\vee$ for the (contravariant) duality functor, and $\Nmix(k)^\oplus$ for the smallest triangulated subcategory of $\Nmot(k)$ which contains $\Nmix(k)$ and is stable under arbitrary direct sums.
\item[(ii)] If $\mathrm{char}(k)=p>0$, then there exists a
  $\bbZ[1/p]$-linear, fully-faithful, symmetric monoidal,
  triangulated functor $\Psi_{1/p}$ making the following diagram commute:
$$
\xymatrix{
\mathrm{Sm}(k) \ar[d]_-{\Sigma^\infty(-_+)_{1/p}} \ar[rrr]^-{X \mapsto \perf_\dg(X)} \ar[drr] &&& \dgcat(k) \ar[d]^-{\mathrm{U}(-)_{1/p}} \\
\mathrm{SH}(k)_{1/p} \ar[d]_-{-\wedge \mathrm{KGL}_{1/p}} && \Nmix(k)_{1/p} \ar[d]_-{(-)^\vee} \ar[r] & \Nmot(k)_{1/p} \ar[d]^-{\uHom(-,\mathrm{U}(k)_{1/p})} \\
\mathrm{Mod}(\mathrm{KGL}_{1/p}) \ar[rr]_-{\Psi_{1/p}} && \Nmix(k)_{1/p}^\oplus \ar[r] & \Nmot(k)_{1/p} \,.
}
$$
\end{itemize}
\end{theorem}
Intuitively speaking, Theorem \ref{thm:bridge1} shows that as soon as we pass to $\mathrm{KGL}$-modules, the commutative world embeds fully-faithfully into the noncommutative world. Consult \cite[\S9.4]{book}, and the references therein, for the construction of the two outer commutative diagrams. The inner commutative squares follow from the combination of Theorem \ref{thm:reduction} with Remark \ref{rk:dualizable}; consult \cite[Thm.~3.1]{Gysin} for details.
\begin{remark}[Morel-Voevodsky's motivic Gysin triangle]\label{rk:motivic1}
Let $X$ be a smooth $k$-scheme, $Z \hookrightarrow X$ a smooth closed subscheme with normal vector bundle $N$, and $j\colon V \hookrightarrow X$ the open complement of $Z$. Making use of homotopy purity, Morel-Voevodsky constructed in \cite[\S3.2]{MV}\cite[\S4]{Voevodsky-ICM} a motivic Gysin triangle
\begin{equation}\label{eq:Gysin-mot2}
\Sigma^\infty(V_+) \stackrel{\Sigma^\infty(j_+)}{\too} \Sigma^\infty(X_+) \too \Sigma^\infty(\mathrm{Th}(N))\stackrel{\partial}{\too} \Sigma^\infty(V_+)[1]
\end{equation}
in $\mathrm{SH}(k)$, where $\mathrm{Th}(N)$ stands for the Thom space of $N$. Since homotopy $K$-theory is an orientable and periodic cohomology theory, $\Sigma^\infty(\mathrm{Th}(N))\wedge \mathrm{KGL}$ is isomorphic to $\Sigma^\infty(Z_+)\wedge \mathrm{KGL}$. Using the commutative diagram \eqref{eq:diagram-1}, we hence observe that the image of \eqref{eq:Gysin-mot2} under the composed functor $\Psi \circ (-\wedge \mathrm{KGL})\colon \mathrm{SH}(k) \to \Nmix(k)^\oplus$ agrees with the dual of the Gysin triangle \eqref{eq:Gysin-mot1}. In other words, the Gysin triangle \eqref{eq:Gysin-mot1} is the dual of the ``$\mathrm{KGL}$-linearization'' of \eqref{eq:Gysin-mot2}.
\end{remark}
\subsection{Voevodsky's motivic category}\label{sub:Voevodsky}
Voevodsky introduced in \cite[\S2]{Voevodsky} the triangulated category of geometric mixed motives $\mathrm{DM}_{\mathrm{gm}}(k)$. By construction, this category comes equipped with a symmetric monoidal functor $M\colon \mathrm{Sm}(k) \to \mathrm{DM}_{\mathrm{gm}}(k)$ and is the natural setting for the study of algebraic cycle (co)homology theories such as higher Chow groups, Suslin homology, motivic cohomology, etc.
\begin{theorem}\label{thm:bridge2}
There exists a $\bbQ$-linear, fully-faithful, symmetric monoidal functor $\Phi_\bbQ$ making~the~following~diagram~commute:
\begin{equation}\label{eq:diagram-3}
\xymatrix{
\mathrm{Sm}(k) \ar[d]_-{M(-)_\bbQ} \ar[rrr]^-{X \mapsto \perf_\dg(X)} \ar[drr] &&& \dgcat(k) \ar[d]^-{\mathrm{U}(-)_\bbQ} \\
\mathrm{DM}_{\mathrm{gm}}(k)_\bbQ \ar[d] && \Nmix(k)_\bbQ \ar[d]_-{(-)^\vee} \ar[r] & \Nmot(k)_\bbQ \ar[d]^-{\uHom(-,\mathrm{U}(k)_\bbQ)} \\
\mathrm{DM}_{\mathrm{gm}}(k)_\bbQ/_{\!\!-\otimes \bbQ(1)[2]} \ar[rr]_-{\Phi_\bbQ} && \Nmix(k)_\bbQ \ar[r] & \Nmot(k)_\bbQ \,.
}
\end{equation}
\end{theorem}
Intuitively speaking, Theorem \ref{thm:bridge2} shows that as soon as we ``$\otimes$-trivialize'' the Tate motive $\bbQ(1)[2]$, the commutative world embedds fully-faithfully into the noncommutative world. Consult \cite[\S9.5]{book}, and the references therein, for the construction of the outer commutative diagram. The inner commutative square follows from the combination of Theorem \ref{thm:reduction} with Remark \ref{rk:dualizable}; consult \cite[Thm.~3.7]{Gysin}.
\begin{remark}[Voevodsky's motivic Gysin triangle]\label{rk:motivic2}
Let $X$ be a smooth $k$-scheme, $Z \hookrightarrow X$ a smooth closed subscheme of codimension $c$, and $j\colon V \hookrightarrow X$ the open complement of $Z$. 
Making use of deformation to the normal cone, Voevodsky constructed in \cite[\S2]{Voevodsky} a motivic Gysin triangle
\begin{equation}\label{eq:Gysin-mot3}
M(V)_\bbQ \stackrel{M(j)_\bbQ}{\too} M(X)_\bbQ \too M(Z)_\bbQ(c)[2c] \stackrel{\partial}{\too} M(V)_\bbQ[1]
\end{equation}
in $\mathrm{DM}_{\mathrm{gm}}(k)_\bbQ$. Using the commutative diagram \eqref{eq:diagram-3}, we observe that the image of \eqref{eq:Gysin-mot3} under the (composed) functor $\Phi_\bbQ \colon \mathrm{DM}_{\mathrm{gm}}(k)_\bbQ \to \Nmix(k)_\bbQ$ agrees with the dual of the rationalized Gysin triangle \eqref{eq:Gysin-mot1}. In other words, the rationalized Gysin triangle \eqref{eq:Gysin-mot1} is the dual of the ``Tate $\otimes$-trivialization'' of \eqref{eq:Gysin-mot3}.
\end{remark}
Let $\mathrm{DM}_{\mathrm{gm}}^{\mathrm{et}}(k)$ be the \'etale variant of $\mathrm{DM}_{\mathrm{gm}}(k)$ introduced by Voevodsky in \cite[\S3.3]{Voevodsky}. As proved in {\em loc. cit.}, $\mathrm{DM}_{\mathrm{gm}}(k)_\bbQ$ is equivalent to $\mathrm{DM}_{\mathrm{gm}}^{\mathrm{et}}(k)_\bbQ$. Consequently, Theorem \ref{thm:bridge2} leads to the following result (see \cite[Thm.~3.13]{Gysin}):
\begin{corollary}[\'Etale descent]
The presheaf of noncommutative mixed motives $\mathrm{Sm}(k)^\op \to \mathrm{NMot}(k)_\bbQ, X \mapsto \mathrm{U}(X)_\bbQ$, satisfies {\em \'etale descent}, \ie for every {\'e}tale cover $\cV=\{V_i \to X\}_{i \in I}$ of $X$, we have an isomorphism $\mathrm{U}(X)_\bbQ \simeq \mathrm{holim}_{n\geq 0} \mathrm{U}(\text{\v{C}}_n\cV)_\bbQ$, where $\text{\v{C}}_\bullet\cV$ stands for the \v{C}ech simplicial $k$-scheme associated to the cover $\cV$.
\end{corollary}
\subsection{Levine's motivic category}\label{sub:Levine}
Levine introduced in \cite[Part I]{Levine} a triangulated category of mixed motives $\cD\cM(k)$ and a (contravariant) symmetric monoidal functor $h\colon \mathrm{Sm}(k) \to \cD\cM(k)$. As proved in \cite[Thm.~4.2]{Ivorra}, the following assignment $h(X)_\bbQ(n) \mapsto \uHom(M(X),\bbQ(n))$ gives rise to an equivalence of categories $\cD\cM(k)_\bbQ \to \mathrm{DM}_{\mathrm{gm}}(k)_\bbQ$ whose precomposition with $h(-)_\bbQ$ is $X \mapsto M(X)^\vee_\bbQ$. Consequently, thanks to Theorem \ref{thm:bridge2}, there exists a $\bbQ$-linear, fully-faithful, symmetric monoidal functor $\Phi_\bbQ$ making the following diagram commute:
\begin{equation}\label{eq:Levine}
\xymatrix{
\mathrm{Sm}(k) \ar[d]_-{h(-)_\bbQ} \ar[rrr]^-{X \mapsto \perf_\dg(X)}&&& \dgcat(k) \ar[dd]^-{\mathrm{U}(-)_\bbQ} \\
\cD\cM(k)_\bbQ \ar[d] &&& \\
\cD\cM(k)_\bbQ/_{\!\!-\otimes \bbQ(1)[2]} \ar[rr]_-{\Phi_\bbQ} && \Nmix(k)_\bbQ \ar[r] & \Nmot(k)_\bbQ \,.\quad \quad 
}
\end{equation}
Note that in contrast with the diagrams of Theorems \ref{thm:bridge1} and \ref{thm:bridge2}, the commutative diagram \eqref{eq:Levine} does {\em not} uses any kind of duality functor.
\subsection{Schur-finiteness conjecture}\label{sec:Schur}
Given a smooth $k$-scheme $X$, the Schur-finiteness conjecture, denoted by $S(X)$, asserts that the mixed motive $M(X)_\bbQ$ is Schur-finite in the sense of Deligne \cite[\S1]{Deligne}. Thanks to the (independent) work of Guletskii \cite{Guletskii} and Mazza \cite{Mazza}, the conjecture $\mathrm{S}(X)$ holds when $\mathrm{dim}(X)\leq 1$, and also for abelian varieties. In addition to these cases, it remains wide open.
\begin{theorem}[{\cite[Thm.~1.1]{Schur}}]\label{thm:fibrations}
Let $q\colon Q \to B$ a flat quadric fibration of relative dimension $d-2$. Assume that $B$ and $Q$ are $k$-smooth and that $q$ has only {\em simple degenerations}, \ie that all the fibers of $q$ have corank $\leq 1$ and that the locus $D \subset B$ of the critical values of $q$ is $k$-smooth. Under these assumptions, the following holds:
\begin{itemize}
\item[(i)] If $d$ is even, then we have $S(B) + S(\widetilde{B}) \Leftrightarrow S(Q)$, where $\widetilde{B}$ stands for the discriminant $2$-fold cover of $B$ (ramified over $D$).
\item[(ii)] If $d$ is odd and $\mathrm{char}(k)\neq 2$, then we have $\{S(V_i)\} + \{S(\widetilde{D}_i)\} \Rightarrow S(Q)$, where $V_i$ is any affine open of $B$ and $\widetilde{D}_i$ is any Galois $2$-fold cover of $D_i:=D\cap V_i$.
\end{itemize}
\end{theorem}
Roughly speaking, Theorem \ref{thm:fibrations} relates the Schur-finiteness conjecture for the total space $Q$ with the Schur-finiteness conjecture for certain coverings/subschemes of the base $B$. Among other ingredients, its proof makes use of Theorem \ref{thm:bridge2} and of the twisted analogue of Theorem \ref{thm:orbifold} (see \S\ref{sec:twisted}). Theorem \ref{thm:fibrations} enables the proof of the Schur-finiteness conjecture in the following new cases:
\begin{corollary}[{\cite[Cor.~1.3 and 1.5]{Fibrations}}]\label{cor:fibrations}
Let $q\colon Q \to B$ be as in Theorem \ref{thm:fibrations}.
\begin{itemize}
\item[(i)] Assume that $B$ is a curve, and that $\mathrm{char}(k)\neq 2$ when $d$ is odd. Under these assumptions, the conjecture $S(Q)$ holds.
\item[(ii)] Assume that $B$ is a surface, that $d$ is odd, that $\mathrm{char}(k)\neq 2$, and that the conjecture $S(B)$ holds. Under these assumptions, the conjecture~$S(Q)$~also~holds.
\end{itemize}
\end{corollary}
Corollary \ref{cor:fibrations}(ii) can be applied, for example, to the case where $B$ is an open subscheme of an abelian surface or smooth projective surface with $p_g=0$ satisfying Bloch's conjecture (see Guletskii-Pedrini \cite[\S4 Thm.~7]{GP}). Recall that Bloch's conjecture holds for surfaces not of general type (see Bloch-Kas-Leiberman \cite{BKL}), for surfaces which are rationally dominated by a product of curves (see Kimura \cite{Kimura}), for Godeaux, Catanese and Barlow surfaces (see Voisin \cite{Voisin2, Voisin}), etc.
\begin{remark}[Bass-finiteness conjecture]
Let $k=\bbF_q$ be a finite field and $X$ a smooth $k$-scheme. The Bass-finiteness conjecture (see \cite[\S9]{Bass}) asserts that the algebraic $K$-theory groups $K_n(X), n \geq 0$, are finitely generated. Thanks to the work of Quillen \cite{Grayson,Quillen2,Quillen1}, the Bass-finiteness conjecture holds when $\mathrm{dim}(X)\leq 1$.

In the same vein, we can consider the {\em mod $2$-torsion} Bass-finiteness conjecture, where $K_n(X)$ is replaced by $K_n(X)_{1/2}$. As proved in \cite{Fibrations}, Theorem \ref{thm:fibrations} and Corollary \ref{cor:fibrations} hold similarly with the Schur-finiteness conjecture replaced by the mod $2$-torsion Bass-finiteness conjecture. As a consequence, we obtain a proof of the (mod $2$-torsion) Bass-finiteness conjecture in new cases.
\end{remark}
\subsection{Rigidity}\label{sub:rigidity2}
Given an integer $n \geq 2$, recall from \cite[\S9.9]{book} the definition of the category of mod-$n$ noncommutative mixed motives $\Nmix(k;\bbZ/n)$. By construction, given smooth proper dg categories $\cA$ and $\cB$, we have isomorphisms
\begin{eqnarray}\label{eq:mod-n}
& \Hom_{\Nmix(k;\bbZ/n)}(\mathrm{U}(\cA), \mathrm{U}(\cB)[-n])\simeq K_n(\cA^\op \otimes \cB; \bbZ/n) & n \in \bbZ\,,
\end{eqnarray}
where the right-hand side stands for mod-$n$ algebraic $K$-theory.
\begin{theorem}[{\cite[Thm.~2.1(ii)]{rigidity}}]\label{thm:rigidity2}
Given an extension of separably closed fields $l/k$, the base-change functor $-\otimes_k l \colon \Nmix(k;\bbZ/n) \to \Nmix(l;\bbZ/n)$ is fully-faithful whenever $n$ is coprime to the characteristic of $k$.
\end{theorem}
Theorem \ref{thm:rigidity2} is the mixed analogue of Theorem \ref{thm:rigidity}. Intuitively speaking, it shows that the theory of mod-$n$ noncommutative mixed motives is ``rigid'' under extensions of separably closed fields. Alternatively, thanks to the isomorphisms \eqref{eq:mod-n}, Theorem \ref{thm:rigidity2} shows that mod-$n$ algebraic $K$-theory is ``rigid'' under extensions of separably closed fields. This is a far-reaching noncommutative generalization of Suslin's celebrated rigidity theorem \cite{Suslin}; consult \cite[\S2]{rigidity} for details and also for applications to equivariant and twisted algebraic $K$-theory. In the particular case of an extension of algebraically closed fields, the commutative counterpart of Theorem \ref{thm:rigidity2} was established by Haesemeyer-Hornbostel in \cite[Thm.~30]{HH}. 
\section{Noncommutative realizations and periods}\label{sec:periods}
In this section we assume that the base field $k$ is perfect. Subsection \S\ref{sec:NCrealizations} is devoted to the noncommutative realizations associated to the (classical) cohomology theories. In \S\ref{sec:periods1}, making use of the noncommutative realization associated to de Rham-Betti cohomology, we extend Grothendieck's theory of periods to the broad noncommutative setting of dg categories. As an application, we prove that (modulo $2\pi i$) Grothendieck's theory of periods is HPD-invariant.
\subsection{Noncommutative realizations}\label{sec:NCrealizations}
Let $F$ be a field of characteristic zero and $(\cC,\otimes, {\bf 1})$ an $F$-linear neutral Tannakian category equipped with a $\otimes$-invertible ``Tate'' object ${\bf 1}(1)$. In what follows, we write $\mathrm{Gal}(\cC)$ for the Tannakian Galois group of $\cC$ and $\mathrm{Gal}_0(\cC)$ for the kernel of the homomorphism $\mathrm{Gal}(\cC) \twoheadrightarrow \bbG_m$, where $\bbG_m$ agrees with the Tannakian Galois group of the smallest Tannakian subcategory of $\cC$ containing ${\bf 1}(1)$. As explained in \cite[\S1-2]{IMRN}, given a cohomology theory $H^\ast\colon \mathrm{Sm}(k) \to \mathrm{Gr}^b_\bbZ(\cC)$, we can consider the associated {\em modified} cohomology theory:
\begin{eqnarray*}
H^{\frac{\ast}{2}}\colon \mathrm{Sm}(k) \too \mathrm{Rep}_{\bbZ/2}(\mathrm{Gal}_0(\cC)) && X\mapsto (\bigoplus_{n\,\,\text{even}} H^n(X), \bigoplus_{n\,\,\text{odd}} H^n(X))
\end{eqnarray*}
with values in the category of finite-dimensional $\bbZ/2$-graded continuous representations of $\mathrm{Gal}_0(\cC)$. Examples of cohomology theories include Nori's cohomology theory $H^\ast_N$ (with values in Nori's Tannakian category of mixed motives \cite[\S8]{Huber}), Jannsen's cohomology theory $H^\ast_J$ (with values in Jannsen's Tannakian category of mixed motives \cite[Part I]{Jannsen1}), de Rham cohomology theory $H^\ast_{dR}$ (with values in the Tannakian category of finite-dimensional $k$-vector spaces), Betti cohomology theory $H^\ast_B$ (with values in the Tannakian category finite-dimensional $\bbQ$-vector spaces), de Rham-Betti cohomology theory $H^\ast_{dRB}$ (with values in the Tannakian category $\mathrm{Vect}(k,\bbQ)$ of triples $(V,W,\omega)$, where $V$ is a finite-dimensional $k$-vector space, $W$ is a finite-dimensional $\bbQ$-vector space, and $\omega$ is an isomorphism $V\otimes_k \bbC\simeq W\otimes_\bbQ \bbC$), \'etale $l$-adic cohomology theory $H^\ast_{l\text{-}\mathrm{adic}}$ (with values in the Tannakian category of finite-dimensional $l$-adic representations of the absolute Galois group of $k$), Hodge cohomology theory $H^\ast_{\mathrm{Hod}}$ (with values in the Tannakian category of mixed $\bbQ$-Hodge structures \cite[\S1]{Steenbrink}), etc; consult \cite[\S2]{IMRN} for details and for further examples. Each one of these cohomology theories gives rise to a modified cohomology theory. 

The (proof of the) next result is contained in \cite[Thms.~1.2, 2.2 and Prop.~3.1]{IMRN}:
\begin{theorem}\label{thm:realizations}
Given a cohomology theory $H^\ast$, there exists an additive invariant
\begin{equation}\label{eq:nc-realization}
H^{\frac{\ast}{2}}_{\mathrm{nc}}\colon \dgcat(k) \too \mathrm{Ind}(\mathrm{Rep}_{\bbZ/2}(\mathrm{Gal}_0(\cC)))\,,
\end{equation}
with values in the category of ind-objects, such that $H^{\frac{\ast}{2}}_{\mathrm{nc}}(\perf_\dg(X))\simeq H^{\frac{\ast}{2}}(X)$ for every smooth $k$-scheme $X$.
\end{theorem}
The additive invariant \eqref{eq:nc-realization} is called the {\em noncommutative realization} associated to the cohomology theory $H^\ast$. Morally speaking, Theorem \ref{thm:realizations} shows that the modified cohomology theories belong not only to the realm of algebraic geometry but also to the broad noncommutative setting of dg categories. This insight goes back to Kontsevich's definition of noncommutative \'etale cohomology theory; see \cite{Kontsevich-talk}. Among other ingredients, the proof of Theorem \ref{thm:realizations} makes use~of~Theorem~\ref{thm:bridge2}.
\begin{remark}[Generalizations]
\begin{itemize}
\item[(i)] In the case where $k$ is of characteristic zero, $\mathrm{Sm}(k)$ can be replaced by the category of $k$-schemes.
\item[(ii)] By construction, \eqref{eq:nc-realization} can be promoted to a localizing invariant. 
\end{itemize}
\end{remark}
The following result describes the behavior of the noncommutative realizations with respect to sheaves of differential operators in characteristic zero\footnote{Consult \cite[Example~2.20]{book} for the description of the behavior of {\em all} additive invariants with respect to sheaves of differential operators in positive characteristic.}:
\begin{theorem}[{\cite[Thm.~3.4]{IMRN}}]\label{thm:D}
Let $k$ be a field of characteristic zero, $X$ a smooth $k$-scheme, and $\cD_X$ the sheaf of differential operators on $X$. Given a cohomology theory $H^\ast$, we have an isomorphism $H^{\frac{\ast}{2}}_{\mathrm{nc}}(\perf_\dg(\cD_X))\simeq H^{\frac{\ast}{2}}(X)$.
\end{theorem}
\begin{example}[Lie algebras]
Let $G$ be a connected semisimple algebraic $\bbC$-group, $B$ a Borel subgroup of $G$, $\mathfrak{g}$ the Lie algebra of $G$, and $U_{\mathrm{ev}}(\mathfrak{g})/I$ the quotient of the universal enveloping algebra of $\mathfrak{g}$ by the kernel of the trivial character. Thanks to Beilinson-Bernstein's celebrated ``localisation'' theorem \cite{BeiBer}, it follows from Theorem \ref{thm:D} that $H^{\frac{\ast}{2}}_{\mathrm{nc}}(U_{\mathrm{ev}}(\mathfrak{g})/I) \simeq H^{\frac{\ast}{2}}_{\mathrm{nc}}(\perf_\dg(\cD_{G/B})) \simeq H^{\frac{\ast}{2}}(G/B)$.
\end{example}
\begin{remark}
Theorem \ref{thm:D} does {\em not} holds for every additive invariant. For example, in the case of Hochschild homology we have $HH_\ast(\perf_\dg(\cD_X))\simeq H_{dR}^{2d-\ast}(X)$ for every smooth affine $k$-scheme $X$ of dimension $d$; see Wodzicki \cite[Thm.~2]{Wodzicki}. Since $H_{dR}^{2d}(X)=0$, this implies that $HH(\perf_\dg(\cD_X))\not\simeq HH(X)$. More generally, we have $HH(\perf_\dg(\cD_X))\not\simeq HH(A)$ for every {\em commutative} $k$-algebra $A$.
\end{remark}
\subsection{Periods}\label{sec:periods1}
In this subsection we assume that the base field is endowed with an embedding $k \hookrightarrow \bbC$. Consider the $\bbZ$-graded $\bbC$-algebra of Laurent polynomials $\bbC[t,t^{-1}]$ with $t$ of degree $1$. Given a triple $(V,W,\omega) \in \mathrm{Vect}(k,\bbQ)$, let us write $\cP(V,W,\omega) \subseteq \bbC$ for the subset of entries of the matrix representations of $\omega$ (with respect to basis of $V$ and $W$). In the same vein, given an object $\{(V_n, W_n, \omega_n)\}_{n \in \bbZ}$ of the category $\mathrm{Gr}_\bbZ^b(\mathrm{Vect}(k,\bbQ))$, let us write $\cP(\{(V_n,W_n,\omega_n)\}_{n \in \bbZ})$ for the $\bbZ$-graded $k$-subalgebra of $\bbC[t,t^{-1}]$ generated in degree $n$ by the elements of the set $P(W_n,W_n, \omega_n)$. In the case of a smooth $k$-scheme $X$, $\cP(X):=\cP(H^\ast_{dRB}(X))$ is called the ($\bbZ$-graded) {\em algebra of periods of $X$}. This algebra, originally introduced by Grothendieck in the sixties, plays a key role in the study of transcendental numbers; consult, for example, the work of Kontsevich-Zagier \cite{KZ}.

Consider the quotient homomorphism $\phi\colon \bbC[t,t^{-1}] \twoheadrightarrow \bbC[t,t^{-1}]/\langle 1 - (2\pi i)t^2\rangle$. The next result extends Grothendieck's theory of periods from schemes to dg categories:
\begin{theorem}[{\cite[Thm.~4.1]{IMRN}}]\label{thm:periods}
There exists an assignment $\cA \mapsto \cP_{\mathrm{nc}}(\cA)$, with $\cP_{\mathrm{nc}}(\cA)$ a $\bbZ/2$-graded $k$-subalgebra of $\bbC[t,t^{-1}]/\langle 1 - (2\pi i)t^2\rangle$, such that $\cP_{\mathrm{nc}}(\perf_\dg(X))$ is isomorphic to $\phi(\cP(X))$ for every smooth $k$-scheme $X$.
\end{theorem}
Morally speaking, Theorem \ref{thm:periods} shows that Grothendieck's theory of periods can be extended from schemes to the broad noncommutative setting of dg categories as long as we work modulo $2\pi i$. Among other ingredients, its proof makes use of the noncommutative realization associated to de Rham-Betti cohomology theory.
\subsubsection{Homological projective duality}
Let $X$ and $Y$ be two HP-dual smooth projective $k$-schemes as in \S\ref{sub:HPD}. Recall from {\em loc. cit.} that the category $\perf(X)$ admits, in particular, a Lefschetz decomposition $\langle \bbA_0, \bbA_1(1), \ldots, \bbA_{i-1}(i-1)\rangle$. Given a linear subspace $L \subset W^\ast$, consider the linear sections $X_L:=X\times_{\bbP(W^\ast)} \bbP(L)$ and $Y_L:=Y \times_{\bbP(W)} \bbP(L^\perp)$. The next result, proved in \cite[Thm.~4.6]{IMRN}, relates the algebra of periods of $X_L$ with the algebra of periods of $Y_L$:
\begin{theorem}[HPD-invariance]\label{thm:HPD2}
Let $X$ and $Y$ be as above. Assume that $X_L$ and $Y_L$ are smooth, that $\mathrm{dim}(X_L)=\mathrm{dim}(X) -\mathrm{dim}(L)$, that $\mathrm{dim}(Y_L)=\mathrm{dim}(Y)- \mathrm{dim}(L^\perp)$, and that the category $\bbA_0$ admits a full exceptional collection. Under these assumptions, the $\bbZ/2$-graded $k$-algebras $\phi(\cP(X_L))$ and $\phi(\cP(Y_L))$ are isomorphic.
\end{theorem}
Roughly speaking, Theorem \ref{thm:HPD2} shows that (modulo $2\pi i$) Grothendieck's theory of periods is invariant under homological projective duality. This result can be applied, for example, to the Veronese-Clifford duality, to the Spinor duality, to the Grassmannian-Pfaffian duality, to the Determinantal duality, etc.

\end{document}

\end{proof}